%% file: article.tex
\documentclass[12pt]{article} 
\usepackage[sectionbib]{natbib}
\usepackage{array,epsfig,fancyhdr,rotating}
\usepackage[]{hyperref}  
\usepackage{sectsty, secdot}
\sectionfont{\fontsize{12}{14pt plus.8pt minus .6pt}\selectfont}
\renewcommand{\theequation}{\thesection\arabic{equation}}
\subsectionfont{\fontsize{12}{14pt plus.8pt minus .6pt}\selectfont}

\usepackage{amsmath}
\usepackage{amssymb}
\usepackage{amsfonts}
\usepackage{multirow}
\usepackage{amsthm}
\usepackage{graphicx}
\usepackage{bm}
\usepackage{subfig}
\usepackage{subfig}
\usepackage{float}
\usepackage{multirow}
\usepackage{hyperref}
\newtheorem{theorem}{Theorem}
\newtheorem{lemma}{Lemma}

\theoremstyle{definition}

\newtheorem{remark}{Remark}

\newcommand\gn{G_{n,q}}
\newcommand\gnqq{G_{n,q-r}}
\newcommand\gp{G_{q}^{(r)}}
\newcommand\Gni{G_{n,q}^{-1}}
\newcommand\gni{G_{n,q}^{-1}}
\newcommand\Hni{H_n^{-1}}
\newcommand\hni{H_n^{-1}}
\newcommand\y{\bm{y}}
\newcommand\z{\bm{\underline{z}}}
\newcommand\ab{\bm{\underline{a}}}
\newcommand\f{\bm{f}}
\newcommand\sfb{\bm{s_f}}
\newcommand\Bt{\bm{B}^T}
\newcommand\bt{\bm{B}^T}
\newcommand\bb{\bm{B}}

\newcommand\bp{B^{(r)}(x)}
\newcommand\bpq{B^{(r)}_q(x)}
\newcommand\bpqt{B^{(r)T}_q(x)}
\newcommand\bqq{B_{q-r}(x)}
\newcommand\bqqt{B_{q-r}^T(x)}

\newcommand\lambdan{\lambda_n}
\newcommand\lambdaP{(\lambdan \Pm)}
\newcommand\lambdap{(\lambdan \Pm)}
\newcommand\cals{\bm{\mathcal{S}}}
\newcommand\Sqt{\cals({q, \underline{\mathbf{t}}})}
\newcommand\sqt{\cals({q, \underline{\mathbf{t}}})}
\newcommand\fph{\hat{f}^{(r)}(x)}
\newcommand\efph{E\left[\hat{f}^{(r)}(x)\right]}
\newcommand\fhp{\hat{f}^{(r)}(x)}
\newcommand\fp{f^{(r)}(x)}
\newcommand\sfp{s_f^{(r)}(x)}
\newcommand\gammab{\underline{\bm{\gamma}}}
\newcommand\gammabt{\underline{\bm{\gamma}}^T}
\newcommand\alphab{\underline{\bm{\alpha}}}
\newcommand\alphabh{\hat{\underline{\bm{\alpha}}}}
\newcommand\alphah{\hat{{\alpha}}}
\newcommand\betab{\underline{\bm{\beta}}}
\newcommand\alphabt{\underline{\bm{\alpha}}^T}
\newcommand\psib{\gammabt\lambdap\hni\gp\hni\lambdap\gammab}
\newcommand\gnh{\gn^{\frac{1}{2}}}
\newcommand\gnmh{\gn^{-\frac{1}{2}}}
\newcommand\tildp{\gnmh\lambdap\gnmh}
\newcommand\Pm{\bm{P_m}}
\newcommand\dr{D^{(r)}}
\newcommand\dtm{\bm{D_m^TD_m}}
\newcommand\dm{\bm{D_m}}
\newcommand\pmh{P_m^{\frac{1}{2}}}
\newcommand\gph{G^{(r)\frac{1}{2}}_q}
\newcommand\calk{\mathcal{K}}
\newcommand\bjqx{B_{j,q}(x)}
\newcommand\mrh{\hat{m}_r(x)}
\newcommand\mh{\hat{m}_0(x)}
\newcommand\hon{h^o_{naive}}
\newcommand\hen{h^e_{naive}}

\begin{document} 

\renewcommand{\baselinestretch}{2}
\newcolumntype{P}[1]{>{\centering\arraybackslash}p{#1}}

\markright{ \hbox{\footnotesize\rm Statistica Sinica
}\hfill\\[-13pt]
\hbox{\footnotesize\rm
}\hfill }

\markboth{\hfill{\footnotesize\rm BRIGHT {ANTWI BOASIAKO} AND JOHN STAUDENMAYER} \hfill}
{\hfill {\footnotesize\rm FILL IN A SHORT RUNNING TITLE} \hfill}

\renewcommand{\thefootnote}{}
$\ $\par


\fontsize{12}{14pt plus.8pt minus .6pt}\selectfont \vspace{0.8pc}
\centerline{\large\bf NAIVE PENALIZED SPLINE ESTIMATORS OF DERIVATIVES }
\vspace{2pt} 
\centerline{\large\bf ACHIEVE OPTIMAL RATES OF CONVERGENCE}
\vspace{.4cm} 
\centerline{Bright Antwi Boasiako and John Staudenmayer} 
\vspace{.4cm} 
\centerline{\it University of Massachusetts, Amherst}
 \vspace{.55cm} \fontsize{9}{11.5pt plus.8pt minus.6pt}\selectfont


\begin{quotation}
\noindent {\it Abstract:}
{\bf This paper studies the asymptotic behavior of penalized spline estimates of derivatives. In particular, we show that simply differentiating the penalized spline estimator of the mean regression function itself to estimate the corresponding derivative achieves the optimal $L_2$ rate of convergence.}

\vspace{9pt}
\par
\end{quotation}\par

\def\thefigure{\arabic{figure}}
\def\thetable{\arabic{table}}

\renewcommand{\theequation}{\thesection.\arabic{equation}}

\fontsize{12}{14pt plus.8pt minus .6pt}\selectfont



\section{Introduction}
\input{background}

\section{Main Results}
\input{main_results}

\section{Simulations}
\label{sec:simulation}
\input{simulations}

\section{Conclusion}
\input{conclusion}

\newpage
\section*{Appendix}
\input{appendix}

\clearpage

\listoffigures
\listoftables





\par


\bibhang=1.7pc
\bibsep=2pt
\fontsize{9}{14pt plus.8pt minus .6pt}\selectfont
\renewcommand\bibname{\large \bf References}
\expandafter\ifx\csname
natexlab\endcsname\relax\def\natexlab#1{#1}\fi
\expandafter\ifx\csname url\endcsname\relax
  \def\url#1{\texttt{#1}}\fi
\expandafter\ifx\csname urlprefix\endcsname\relax\def\urlprefix{URL}\fi

\bibliography{report}
\bibliographystyle{apalike}

\vskip .65cm
\noindent
Bright Antwi Boasiako
\vskip 2pt
\noindent
E-mail: (bantwiboasia@umass.edu)
\vskip 2pt

\noindent
John Staudenmayer
\vskip 2pt
\noindent
E-mail: (jstauden@umass.edu)

\end{document}

%% file: background.tex
We consider the situation where data $\{x_i, y_i\}_{i=1}^n$, sampled from the model:
\begin{equation}
    \label{eqn:1}
    y_i = f(x_i)+\varepsilon_i,\ \ \forall i=1, 2, \dots, n
\end{equation}
with $f \in \mathcal{C}^p\mathcal{(K)}$, the space of functions with $p$ continuous derivatives over $\mathcal{K}=\left[0, 1\right]$, the $x_i$'s, for $x_i\in \mathcal{K}$, are either random or deterministic, and $\varepsilon_i$'s are independent and identically distributed random error terms with $E[\varepsilon_i]=0$ and $Var[\varepsilon_i] = \sigma^2$. There are many cases when it is of interest to estimate some derivative of the mean regression function $f$, with minimal assumptions on the functional form of $f$. For example, in human growth studies, the first derivative of the function relating height and age indicates the speed of growth (\citealp[]{muller1988, ramsay2002}). Additionally, \cite{Chaudhuri1999} apply derivative estimation in the development of a visual mechanism for studying curve structures, and \cite{Park2008} compare regression curves using those structures.
In economics, derivatives are used to calculate the marginal propensity to consume, which measures the effect of changes in disposable income on personal consumption (\citealp[]{FISHER2020103218}). In addition, average derivatives of mean regression functions are used to empirically validate the so-called ``law of demand" (\citealp[]{hardle1989}). It is sufficient for a random matrix composed of average derivatives to be positive definite for the law of demand to hold (\citealp[]{HILDENBRAND1989251}).
In nonparametric regression itself, estimates of derivatives of the true function, $f$, are used in plug-in bandwidth selection techniques such as in local polynomial regression (\citealp[]{ruppert_plugin_1995}) and to construct confidence bands for nonparametric estimators (\citealp[]{Eubank1993}).

Previous work has taken three major approaches to estimate derivatives of functions nonparametrically: local polynomial regression, empirical derivatives, and spline-based methods. In local polynomial regression, a derivative of $f(x)$ can be estimated using a coefficient of the fitted local polynomial at $x$ \cite[p~22]{fan1996local}. Empirical methods generally transform the data and smooth difference-based estimates of derivatives. Earlier works include \cite{muller1987}, which used Kernel-based approaches to estimate the derivatives of the mean regression function while utilizing difference quotients to identify the best kernel bandwidth via cross validation. More recently, \cite{brabanter2013} used symmetric difference quotients to estimate derivatives of mean regression functions and showed that their approach improved the asymptotic order of the variance. Spline-based methods use the fact that splines are piecewise polynomials. As a result, differentiating the basis function with respect to the covariate gives a basis function for the derivative, and an estimate of that function then can be obtained using a subset of the estimated coefficients (\citealt[p~115]{deboor1978}, \citealp[]{eilers1996, eilers2010}).

While it is straightforward to compute nonparametric derivative estimates, the challenge is that those estimates also require some sort of regularization to balance estimation bias and overfitting, and methods to choose that regularization are usually designed for estimating the function itself, not derivatives (\citealp[]{ruppert_wand_carroll_2003, eilers1996}). Several authors have designed methodology to address that problem (e.g. \citealp[]{charnigo2011, simpkin2013}), but it has also been suggested that methods that choose the amount of smoothing for the derivatives as if the function were of interest often work well in practice (\citealp[p~154]{ruppert_wand_carroll_2003}, \citealp[]{craven_smoothing_1978}). We call such methods \textit{naive}. In this paper, we add evidence to that debate by exploring the asymptotic behavior of naive nonparametric derivative estimators, focusing on penalized splines.

The past few decades have seen some progress in related areas. In local polynomial regression, results from \cite{ruppert_wand_1994} can be used to show that when a bandwidth is chosen to minimize the integrated mean squared error (IMSE) of a $p^{th}$ degree polynomial estimate of $f$, and that bandwidth is used to estimate $f^{(r)}$ ($p \ge r$, i.e. a naive estimator), then the derivative estimate’s IMSE converges at an optimal rate if $r$ is even. Otherwise, the naive bandwidth over- or under-smooths. In particular, a naive estimator of the first derivative using cubic local polynomials under-smooths. A derivation of this is given in the appendix. Also, \cite{brabanter2013} showed that the asymptotic order of the bias of empirical estimators does not depend on the order of the derivative being estimated, and they employ a method by \cite{debrabanter2011} to deal with correlations that result from creating the empirical dataset for the derivative. \cite{daitong2016} generalized those results by considering linear combinations of observations to better fit both interior and boundary points. They also demonstrated that their method achieves optimal rates of convergence (\citealp[]{stone1982}). The asymptotic properties of two types of spline-based derivative estimators have been considered too. \cite{wangwahba1990} studied smoothing splines and found that the optimal smoothing parameter depends on the order of the derivative being estimated. In contrast, \cite{zhou2000yhh} considered the asymptotics of regression spline based derivative estimators where the number of knots increases with the sample size. They showed that the MSE goes to zero at the optimal rate (\citealp[]{stone1982}), and the required rate of increase in the number of knots does not depend on the order of the derivative.

\begin{sloppypar}
Somewhat surprisingly though, comparable results about penalized spline estimators of derivatives do not seem to exist. Building on work on the asymptotics of penalized spline estimators of functions (\citealp[]{xiao2019}), we derive the apparently new result that naive methods to estimate derivatives with penalized splines achieve optimal global $L_2$ rates of convergence (\citealp[]{stone1982}). 
\end{sloppypar}

The rest of the paper is structured as follows: in Section 2, we give some background on splines, penalized splines and the naive derivative estimator. In Section 3, we present our main result and remark that depending on the rate at which the number of knots increases with the sample size, $n$, the $L_2$ convergence of the naive estimator is similar to that of regression splines or smoothing splines.  In Section 4, we present a simulation study of the $L_2$ rate of convergence of the naive estimator, and we conclude with a discussion in Section 5.  Proof of our main result and other technical lemmas are given in the appendix.

\section{Penalized Splines \& the Naive Derivative Estimator}
\subsection{Splines}
\begin{sloppypar}
Splines provide a flexible mechanism to estimate derivatives of the the mean regression function $f$, and in the case of estimating the function itself, they have been shown to do so at the best possible rates of convergence (\citealp[]{xiao2019, zhou1998, stone1982}). 
A spline is a piece-wise polynomial with continuity conditions at the points where the pieces join together (called knots).\\
More specifically, for $q\ge2$, we let $$\sqt = \left\{ s\in \mathcal{C}^{q-2}(\calk): s\text{ is a $q$-order polynomial on each }[t_i, t_{i+1}] \right\}$$ be a space of $q-$order splines over $\calk = [0, 1]$ with knot locations $\underline{\mathbf{t}} = (t_0, t_1, \dots, t_{K+1})$ where $t_0 = 0$, $t_{K+1}=1$ and $t_i < t_j$ $\forall_{i < j}$. For $q=1$, $\sqt$ consists of step functions with jumps at the knots.
\end{sloppypar}
This space has a number of equivalent bases and one notable for having stable numerical properties is the \textit{B-Spline} basis (\citealp[]{deboor1978, ruppert_wand_carroll_2003, schumaker_2007}).

\cite{deboor1978} defines the $q^{th}$ order B-spline basis function $\bjqx$ over the knot locations $\underline{\mathbf{t}}$ through a recurrence relation. \cite{eilers2010} show that when the distance between the knots is constant, $\bjqx$ reduces to 
$$\bjqx = \frac{(-1)^q\Delta^q(x-t_j)_+^{q-1}}{(q-1)!h^{q-1}}$$
where $\Delta$ is the backward difference operator ($\Delta t_j=t_j - t_{j-1}$), and $h$ is the common distance between the knots. Observe that $\bjqx$, in this case, is a rescaled $q$-order difference of truncated polynomials.
To get a complete set of B-Spline basis, we need $2q$ extra knots with $q$ knots on each side of $\left[0, 1\right]$. This is referred to as the expanded basis (\citealp[]{eilers2010}).

Without losing generality, we will assume a B-Spline basis for $\sqt$ for the rest of this paper. We refer the reader to \cite{deboor1978, schumaker_2007, eilers1996} for an introduction to B-Splines, and \cite{eilers2010} for how the B-Spline basis compares to the Truncated Polynomial Functions (TPF) on metrics including fit quality, numerical stability, and multidimensional smoothing.

\subsection{Penalized Splines \& the Naive Estimator}
Penalized splines are often viewed as a compromise between regression and smoothing splines because they combine penalization and low rank bases to achieve computational efficiency. They vary slightly based on the basis functions used and the object of penalization. For example, P-Splines (\citealp[]{eilers1996}) use B-Spline basis functions and penalize differences of the coefficients to a specific order. In this section, we will focus on P-Splines. Our later results hold for the general penalized spline estimator defined by \cite{xiao2019}.

A P-Spline estimator of $f$ in \eqref{eqn:1} based on an iid sample of size $n$ finds a \textit{spline} function $g(x) = B(x)\alphab$, that minimizes: 
\begin{equation}
    \label{eqn:obj}
    Q(\alphab, \lambdan) = \frac{1}{n} \sum_{i=1}^n \left(y_i - B(x_i)\alphab\right)^2 + \lambdan \alphabt\Pm\alphab
\end{equation}
where, $\bm{\alphab}=\left(\alpha_1, \alpha_2, \dots, \alpha_{K+q}\right)$ is a vector of coefficients, and $B(x_i) = \left[B_{1,q}(x_i), B_{2,q}(x_i), \dots, B_{K+q,q}(x_i)\right] \in \mathbb{R}^{K+q}$ is a vector of basis functions at $x_i$, for $i = 1, 2, \dots, n$. The penalty matrix $\Pm=\dtm \in \mathbb{R}^{(K + q) \times (K + q)}$ where $\dm\alphab = \Delta^m\alphab$ is a vector of $m^{th}$ order differences of $\alphab$. Finally, $\lambdan \ge 0$ is the smoothing parameter and needs to be chosen. Three prevalent methods for choosing $\lambdan$ are Generalized Cross Validation (GCV), Maximum Likelihood (ML), and Restricted (or residual) Maximum Likelihood (REML); we refer the reader to \cite{wood_2017} chapter 4 and \cite{ruppert_wand_carroll_2003} Chapters 4 and 5 for details.
 
Minimizing \eqref{eqn:obj} with respect to $\alphab$ gives
$\alphabh=\left(\frac{\bt\bb}{n} + \lambdan \Pm\right)^{-1}\frac{\bt y}{n}$ which results in $\hat{f}(x) = B(x)\alphabh$. Here, $\bm B=\left[B(x_1), B(x_2), \dots, B(x_n)\right]^T \in \mathbb{R}^{n\times K+q}$.
From $\hat{f}$, we can derive the naive estimator of the $r^{th}$ derivative of $f$ as follows:
\begin{eqnarray}
\label{eqn:naive}
    \hat{f}^{(r)}(x) &=& \frac{d^{(r)}}{dx}B(x)\alphabh \nonumber \\
    &=& \frac{d^{(r)}}{dx}\left(\sum_{j=1}^{K+q}\alphah_jB_{j, q}(x)\right) \nonumber \\
    &=& \displaystyle\sum_{j=1}^{K+q-r}\alphah_j^{(r)}B_{j, q-r}(x)
\end{eqnarray}

where $\alphah_j^{(r)} = (q-r)\frac{\left(\alphah_{j+1}^{(r-1)}-\alphah_{j}^{(r-1)}\right)}{t_j-t_{j-q+r}}$, with $\alphah_j^{(0)} = \alphah_j$ for $1\le j\le K+q-r$, and $r = 1, 2, \dots, q - 2$. (\citealp[]{deboor1978, zhou2000yhh}).

\cite{xiao2019} showed that under some conditions on the distribution of the knots and $\lambdan$, $\hat{f}$ achieves the optimal $L_2$ rate of convergence (\citealp[]{stone1982}) to the true $f$ but they do not discuss derivative estimators. We extend this result to show that under same conditions that do not depend of the order of the derivative, the naive derivative estimator $\hat{f}^{(r)}$ of $f^{(r)}$ achieves optimal $L_2$ rates of convergence.

%% file: main_results.tex
In this section, we provide our main result in Theorem \ref{thm:1} and remark on how this result relates to regression and smoothing splines. Note that the findings in this section apply to the general penalized spline estimator as defined by \cite{xiao2019}. This general estimator is based on the realization that the various types of penalized splines differ mainly by their penalty matrices. However, the eigenvalues of the penalty matrices decay at similar rates, making their unified asymptotic study tractable. We refer the interested reader to a derivation of the decay rates of various penalty matrices in \cite{xiao2019}.
\subsection{Notation}
We start by defining the following notations relating to norms and limits. For a real matrix $A$, $||A||_\infty = \displaystyle \max_i \sum_{j} |a_{ij}|$ is the largest row absolute sum. $||A||_2$ is the operator norm of $A$ induced by the vector norm $||.||_2$. $||A||_F = \sqrt{tr(A^TA)}$ is the Frobenius norm. For a real vector, $||\ab|| = \displaystyle \max_i |a_i|$. For a real-valued function $g(x)$ defined on $\calk\subset \mathbb{R}$, $||g|| = \displaystyle \sup_{x\in\calk} |g(x)|$ and $||g||_{L2}=\left(\displaystyle\int_{x\in\mathcal{K}}\left(g(x)\right)^2dx\right)^{1/2}$ is the $L_2$-norm of $g$. For two real sequences $\{a_n\}_{n\ge1}$ and $\{b_n\}_{n\ge1}$, $a_n\sim b_n$ means $\displaystyle\lim_{n\to\infty}\frac{a_n}{b_n} = 1$.

\subsection{Assumptions}
Next, we state assumptions on the knot placement and penalty matrix. We note that these assumptions are the same as those made in \cite{xiao2019} for the asymptotic analysis of estimates of functions rather than derivatives.

\begin{enumerate}
    \item $K=o(n)$.
    \label{assumption:1}
    \vspace{15pt}
    \item $\displaystyle \max_{1\le i\le K} |h_{i+1}-h_i| = o(K^{-1})$, where $h_i = t_i - t_{i-1}$.
    \label{assumption:2}
    \vspace{15pt}
    \item $\frac{h}{\displaystyle \min_{1\le i\le K}h_i} \le M$, where $h = \displaystyle \max_{1\le i\le K}h_i$ and $M>0$ is some predetermined constant.
    \label{assumption:3}
    \vspace{10pt}
    \item For a deterministic design,
    $$\displaystyle \sup_{x\in [0, 1]}|Q_n(x) - Q(x)| = o(K^{-1})$$ 
    where $Q_n(x)$ is the empirical CDF of $x$ and $Q(x)$ is a distribution with continuously differentiable positive density $q(x)$.
    \label{assumption:4}
    \item The penalty matrix $\Pm$ is a banded symmetric positive semi-definite square matrix with a finite bandwidth and $||\Pm||_2 = O(h^{1-2m})$. 
    This assumption is similar to Assumption 3 of \cite{xiao2019} where it is stated in terms of the eigenvalues of $\Pm$. This assumption is verifiable for P-Splines, O-Splines and T-Splines. See Propositions 4.1 and 4.2 of \cite{xiao2019}.
    \label{assumption:5}
    \item $\lambdan = o(1).$
    \label{assumption:6}
\end{enumerate}

Assumptions \eqref{assumption:2} and \eqref{assumption:3} are necessary conditions on the placements of the knots and also imply that $h\sim K^{-1}$. This ensures that $M^{-1} < Kh < M$ and is necessary for numerical computations (\citealp[]{zhou1998}).

\begin{theorem}
\label{thm:1}
    Let the mean regression function in \eqref{eqn:1} be such that $f \in \mathcal{C}^p\mathcal{(K)}$. Under Assumptions \eqref{assumption:1} - \eqref{assumption:6} above, and for $m\le \min(p, q)$:
    \begin{eqnarray*}
        \mathbb{E}\left(|| \hat{f}^{(r)} - f^{(r)} ||_{L_2}^2\right) &=&
        O\left(\frac{K_e}{n}\right) + 
        O\left(K^{-2(q-r)}\right) + o(K^{-2(p-r)}) \\
        && + O\{ \min(\lambdan^2K^{2m+2r}, \lambdan K^{2r}) \}
    \end{eqnarray*}
    where $K_e = \min\left\{K^{2r+1}, K^{2r}\lambdan^{-1/2m}\right\}$ and $r = 1, 2, \dots, q-2$.
\end{theorem}

The proof of the theorem is given in the appendix.

\subsection{\textbf{Remarks}}

\begin{remark}
\label{rmk:1}
\noindent The asymptotics of penalized splines are either similar to those of regression splines or smoothing splines depending on how fast the number of knots increases as the sample size increases (\citealp[]{Claeskens2009, xiao2019}). This creates two scenarios: the small number of knots scenario with asymptotics similar to regression splines and the large number of knots scenario with asymptotics similar to smoothing splines. We explore the rates of convergence of the naive estimator under each of these scenarios in Remarks 1a and 1b below. 

\noindent {\bf Remark 1a} (Small number of knots scenario):
Suppose the mean regression function is $q$-times continuously differentiable, where $q$ is the order of the spline used to estimate $f$. Thus, $f\in \mathcal{C}^q(\mathcal{K})$. Also suppose $\lambdan K^{2m} = O(1)$, then
\begin{eqnarray*}
\mathbb{E}\left(|| \hat{f}^{(r)} - f^{(r)} ||_{L_2}^2\right) &=& 
    O\left(\frac{K_e}{n}\right) +
O\left(K^{-2(q-r)}\right) + o(K^{-2(p-r)}) \\
        && + O\{ \min(\lambdan^2K^{2m+2r}, \lambdan K^{2r}) \}\\
        &=& O\left(\frac{K^{2r+1}}{n}\right) + O\left(K^{-2(q-r)}\right)
        + O( \lambdan^2K^{2m+2r}).
    \end{eqnarray*}
Choosing $K$ such that $K\sim n^{\frac{1}{2q+1}}$ and $\lambdan = O(n^{-(q+m)/(2q+1)})$, the estimator
$\hat{f}^{(r)}$ of $f^{(r)}$ converges at the optimal $L_2$ rate of $n^{-\frac{(q-r)}{2q+1}}$. In the above, we have used the fact that $p=q$ and that $\min\left\{\lambdan^2K^{2m+2r}, \lambdan K^{2r}\right\} = \lambdan^2K^{2m+2r}$, $K_e = K^{2r+1}$ for $\lambdan K^{2m} = O(1)$. We note that the $\lambda_n$'s rate of decrease does not depend on $(r),$ the order of the derivative.

\noindent {\bf Remark 1b} (Large number of knots scenario):
Suppose $f\in \mathcal{C}^m(\mathcal{K})$, and there exists a sufficiently large constant, $C$, independent of $K$ such that for $K\ge C^{1/2m}\lambdan^{-1/2m} = C^{1/2m}n^{\frac{1}{2m+1}}$, with $m \le q$, we have
\begin{eqnarray*}
        \mathbb{E}\left(|| \hat{f}^{(r)} - f^{(r)} ||_{L_2}^2\right) &=& 
        O\left(\frac{K_e}{n}\right) + 
        O\left(K^{-2(q-r)}\right) + o(K^{-2(p-r)}) \\
        &&+ O\{ \min(\lambdan^2K^{2m+2r}, \lambdan K^{2r}) \}\\
        &=& O\left(\frac{K^{2r}\lambdan^{-1/2m}}{n}\right) + O\left(K^{-2(q-r)}\right) 
        + o\left(K^{-2(m-r)}\right) \\
        && + O( \lambdan K^{2r}).
    \end{eqnarray*}
Choosing $\lambdan$ such that $\lambdan\sim n^{-2m/(2m+1)}$, the estimator $\hat{f}^{(r)}$ of $f^{(r)}$ converges at the optimal $L_2$ rate of $n^{-\frac{(m-r)}{2m+1}}$. Again, we note that the $\lambda_n$'s rate of decrease does not depend on $(r).$

\end{remark}

\begin{remark}
\label{rmk:2}
\noindent While the naive estimator of the derivative achieves an optimal rate of convergence, that does not mean that the naive approach is optimal in a finite sample. We compare the performance of the naive estimator to an ``oracle estimator" that minimizes mean integrated squared error in Section \hyperref[sec:finitesample]{4.1.4}.

\end{remark}

\begin{remark}
\label{rmk:3}
\noindent The theorem is derived under conditions on the growth in the number of knots, the spacings between them, and the smoothing parameter ($\lambda_n$). Specific rates of growth for $K$ and for $\lambda_n$ in Remarks 1a and 1b led to optimal rates of convergence. That said, it is not clear whether standard ways of choosing smoothing parameters would lead to optimal rates of converged. This too is explored in Section \hyperref[sec:simulation]{4}.
\end{remark}

%% file: simulations.tex
\subsection{Overview}

In this section, we present a simulation to assess the naive estimator's rate of convergence and its finite-sample performance. The simulation is divided into three parts. The first part examines the $L_2$ rates of convergence of the naive estimator when GCV and REML are used to choose the smoothing parameter. The second part of this section focuses on the finite sample performance of the naive estimator. We compared it to an ``oracle'' method that uses knowledge of the true function (or derivatives) to choose the optimal smoothing parameter. That ``oracle'' method is not a practical estimator, but it provides an upper bound benchmark for P-spline performance. Finally, the third part of this section compares the naive method to other derivative estimation methods in the literature. 

\begin{figure}[H]
    \centering
    \subfloat[\centering {\scriptsize Mean regression function $f(x)$}]{{\includegraphics[width=4.6cm]{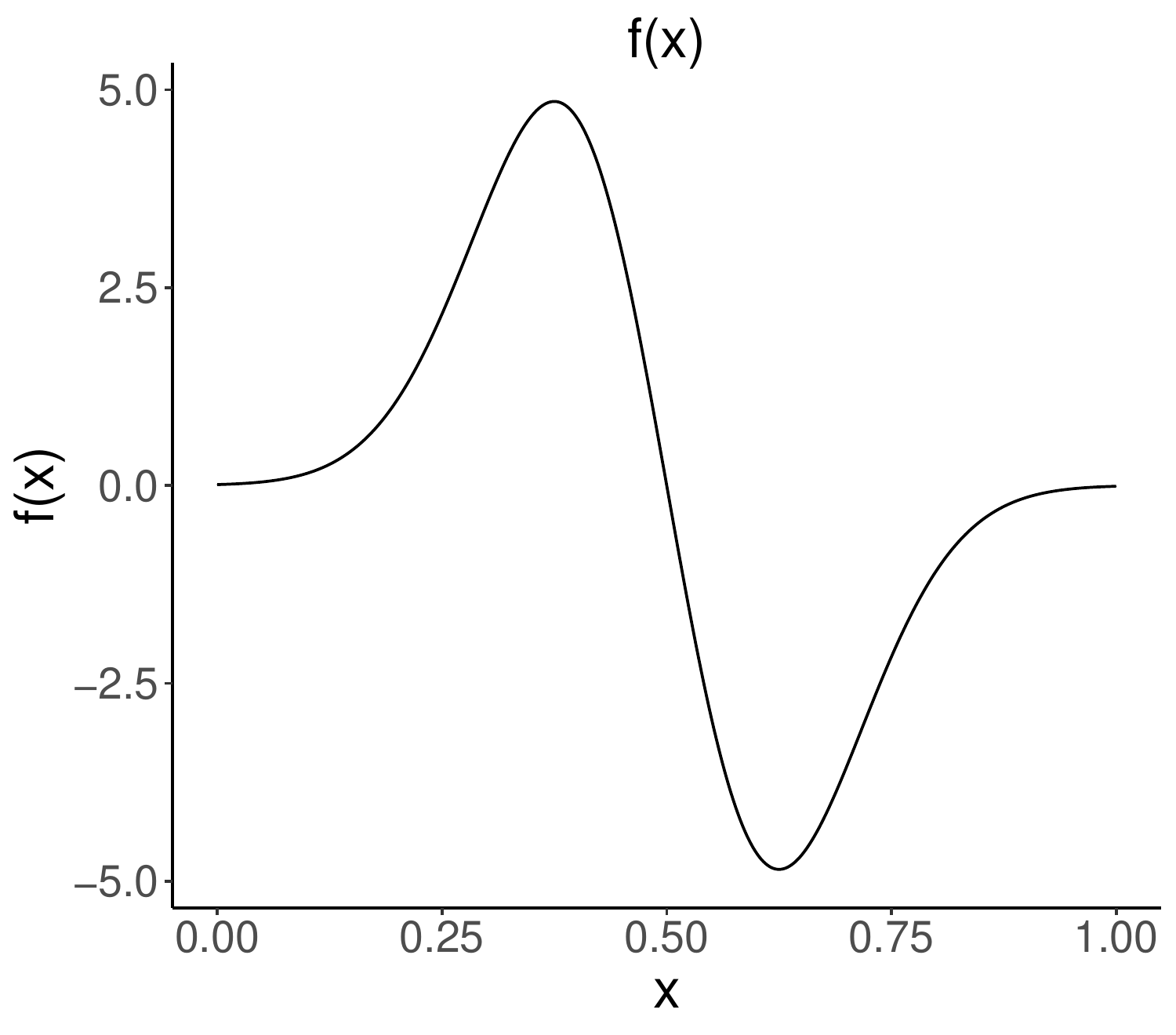} }}%
    \subfloat[\centering \scriptsize First derivative of $f(x)$]{{\includegraphics[width=4.6cm]{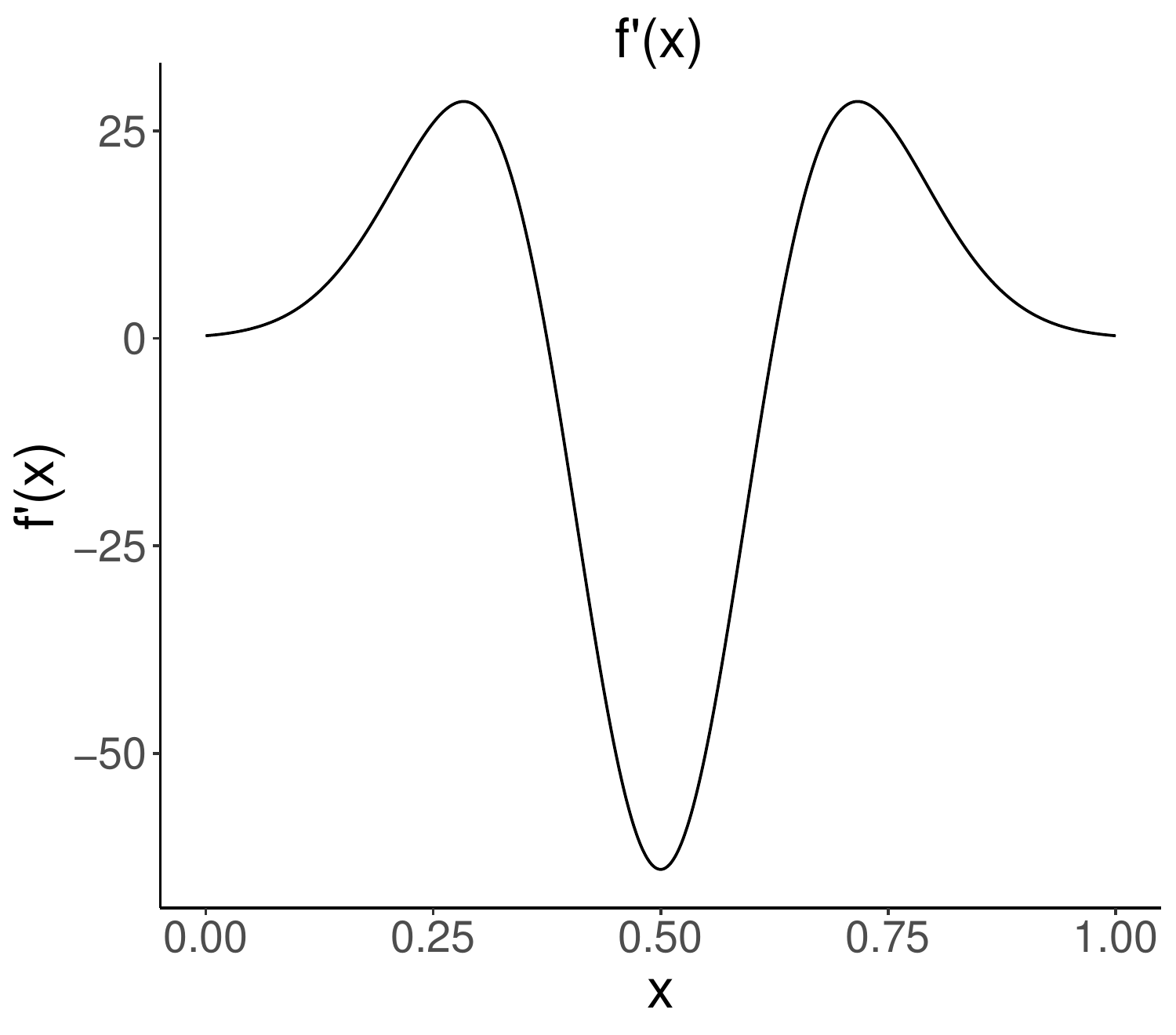} }}%
    \subfloat[\centering \scriptsize Second derivative of $f(x)$]{{\includegraphics[width=4.6cm]{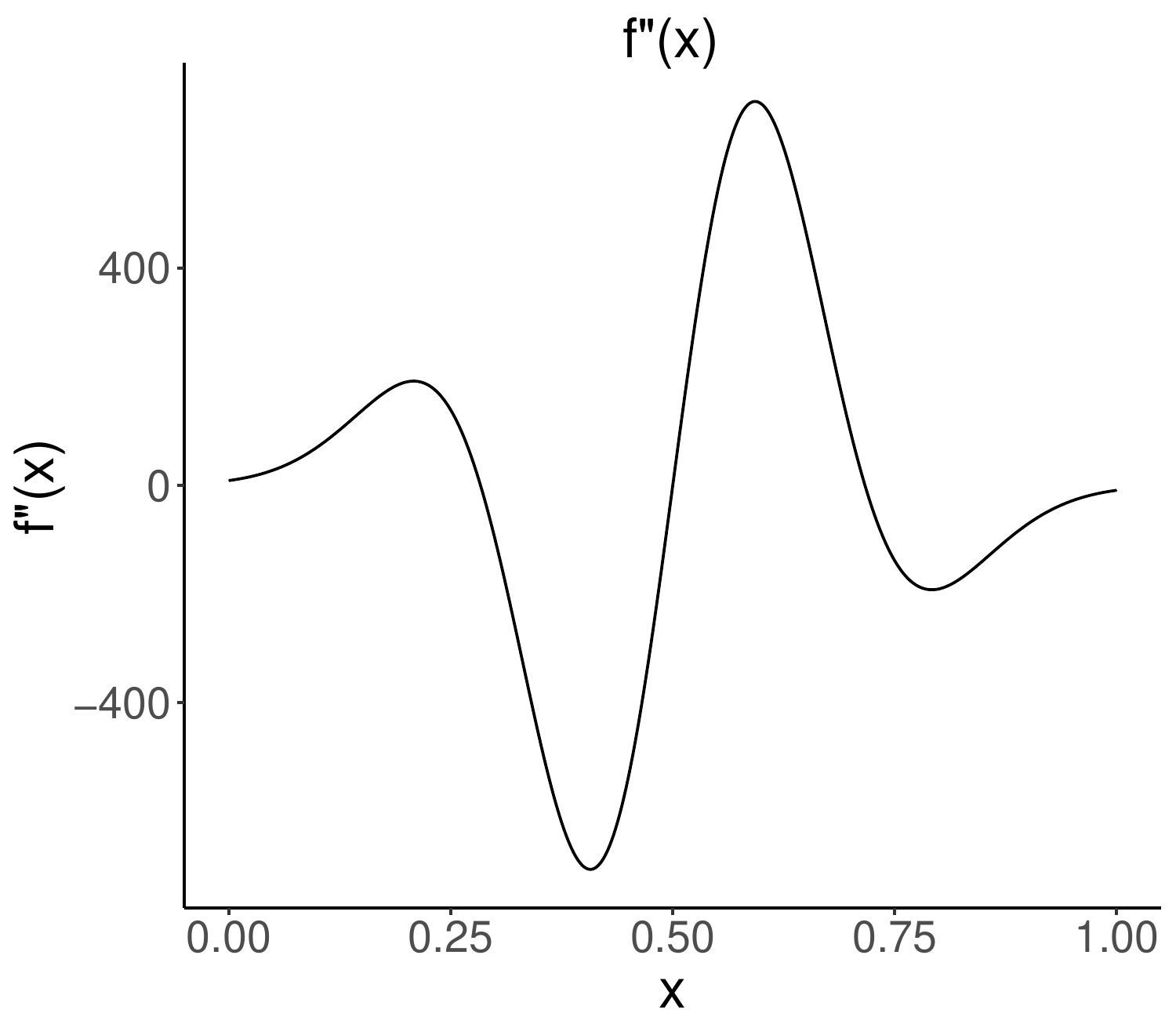} }}%
    \caption{Mean regression function with its first two derivatives.}%
    \label{fig:funcs}%
\end{figure}

Except where noted, we use the same mean regression function $f$ as \cite{brabanter2013}. We simulated data $\left\{x_i, y_i\right\}_{i=1}^n$ from the model:
$$Y_i = f(x_i) + \varepsilon_i, \ \forall\ 1\le i\le n$$ 
where $x_i$'s are a grid over $\mathcal{K} = \left[0,1\right]$, $\varepsilon_i$'s are iid with $\varepsilon_i \sim N(0, \sigma^2 = 0.1^2)$ and
\begin{equation}
\label{eqn:barbranter}
    f(x)= 32e^{-8(1-2x)^2}\left(1-2x\right)
\end{equation}
Figure \ref{fig:funcs} shows the mean regression function in \eqref{eqn:barbranter} and its first two derivatives. We use a range of sample sizes as shown in the results.

As discussed in \cite{xiao2019}, \cite{Claeskens2009}, and our Remark \eqref{rmk:1}, the asymptotics of the penalized spline estimator are similar to those of Regression Splines (small K scenario) or Smoothing Splines (large K scenario), depending on the rate at which the number of knots, $K$, increases with the sample size, $n$. In our simulation, we considered these two scenarios: when $K$ increases slowly with $n$, and when $K$ increases at a faster rate with $n$. For the slow $K$ scenario we use $K\sim n^{\frac{1}{2p+1}},$ and $K$ in the fast scenario is chosen such that $K\ge C^{1/p}\lambdan^{-1/2p}$ for some large constant, $C$.

We investigated the $L_2$ rate of convergence for the first two derivatives of the mean regression function in \eqref{eqn:barbranter} using a P-Spline with $2^{nd}$ ($m=2$) order penalty (\citealp[]{eilers1996}). Note that with $m=2$, the equivalent kernel methodology (\citealp[]{silverman1984}, Lemma 9.13 of \citealp[]{xiao2012local}) implies that the assumed differentiability of $f$ is $p=2m = 4$.

\cite{stone1982} provided optimal rates of convergence for non-parametric regression estimators. The optimal rate of convergence for a non-parametric estimator of the $r^{th}$ derivative of $g:\mathbb{R}^d \rightarrow \mathbb{R}$ where $g \in \mathcal{C}^p$ is given by $n^{-\frac{p-r}{2p+d}}$, in our simulations, we have the optimal $L_2$ rate of convergence for estimating the $r^{th}$ derivative of $f$ as: 
$$ n^{-\frac{p-r}{2p+d}} = n ^ {-\frac{4-r}{2 \times 4 + 1}} = n^{-\frac{1}{9}(4-r)}$$

\subsubsection{$L_2$ Convergence of the Naive Estimator}

\begin{figure}[H]
    \centering
    \subfloat[\centering \scriptsize $L_2$ rate of convergence for $\hat{f}$]{{\includegraphics[width=4.6cm]{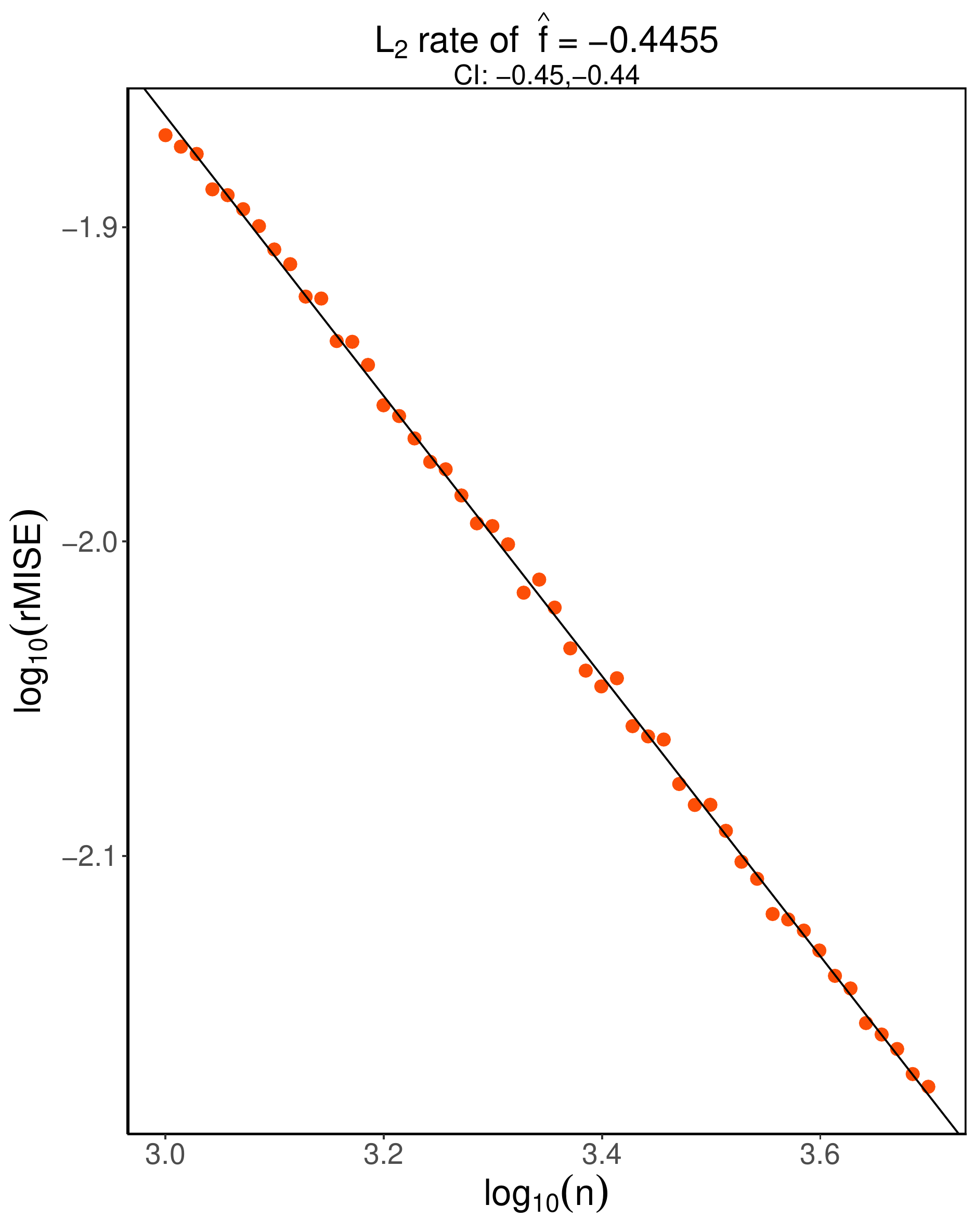} }}%
    \subfloat[\centering \scriptsize $L_2$ rate of convergence for $\hat{f}'$]{{\includegraphics[width=4.6cm]{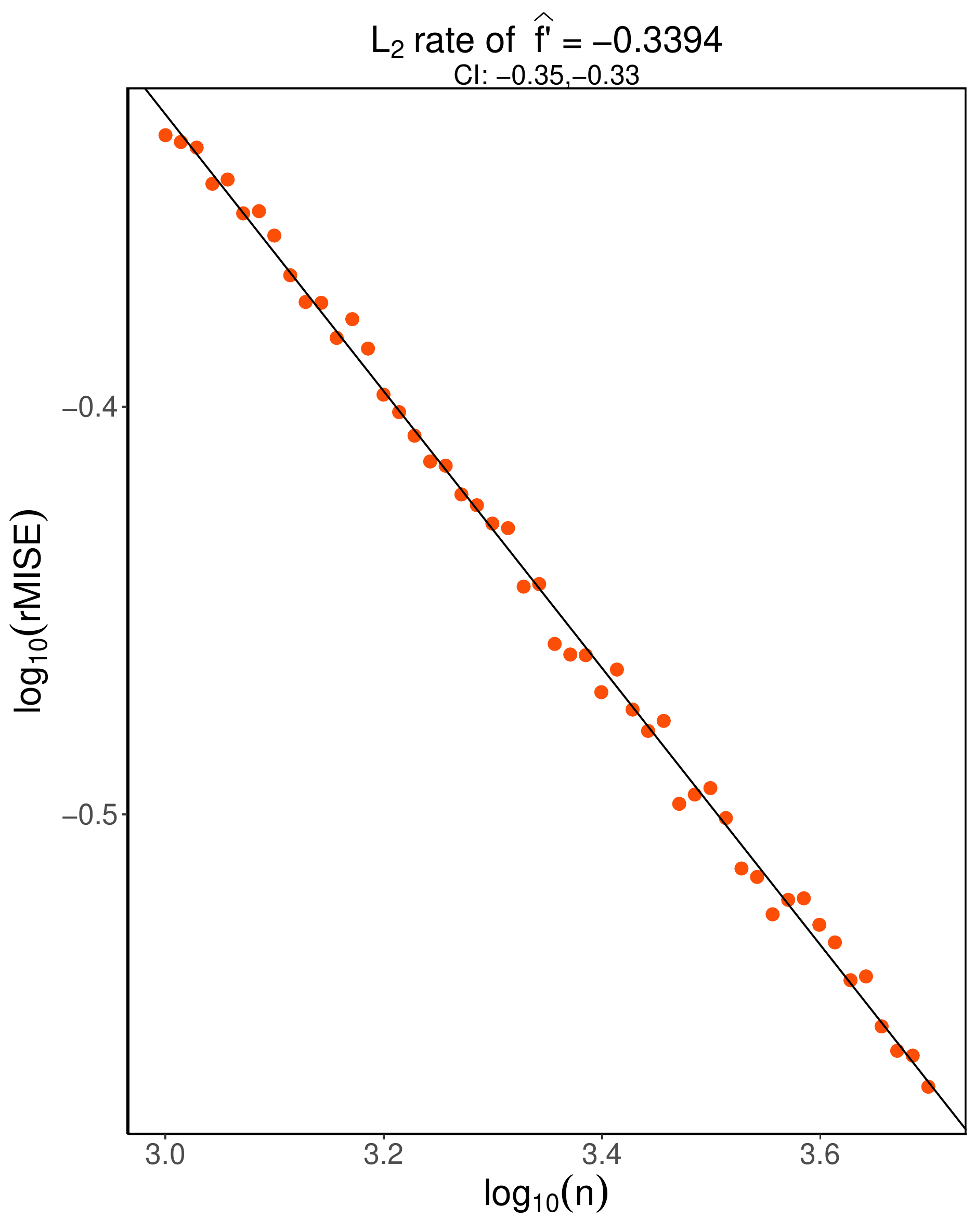} }}%
    \subfloat[\centering \scriptsize $L_2$ rate of convergence for $\hat{f}''$]{{\includegraphics[width=4.6cm]{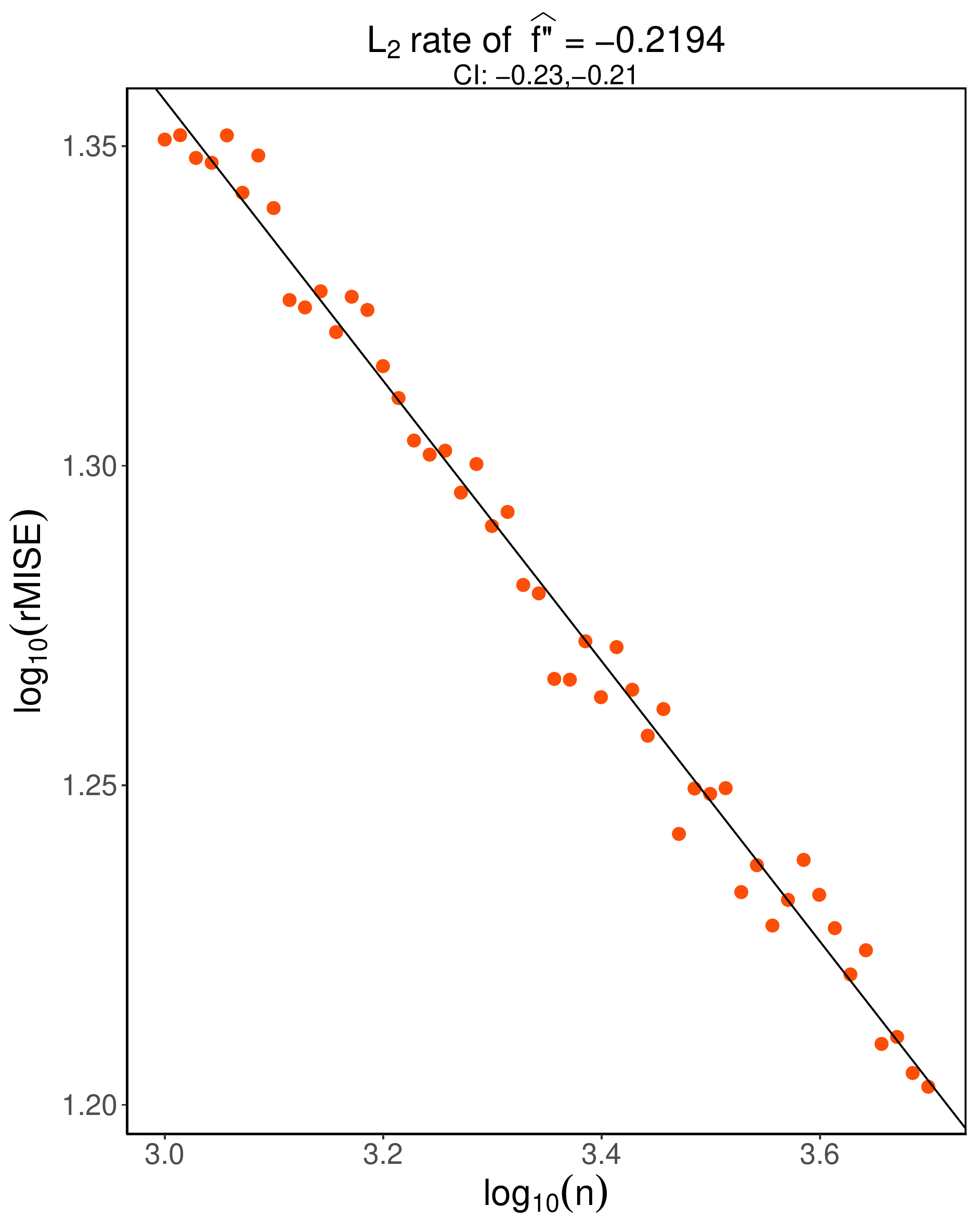} }}%
    
    \centering
    \subfloat[\centering \scriptsize $L_2$ rate of convergence for $\hat{f}$]{{\includegraphics[width=4.6cm]{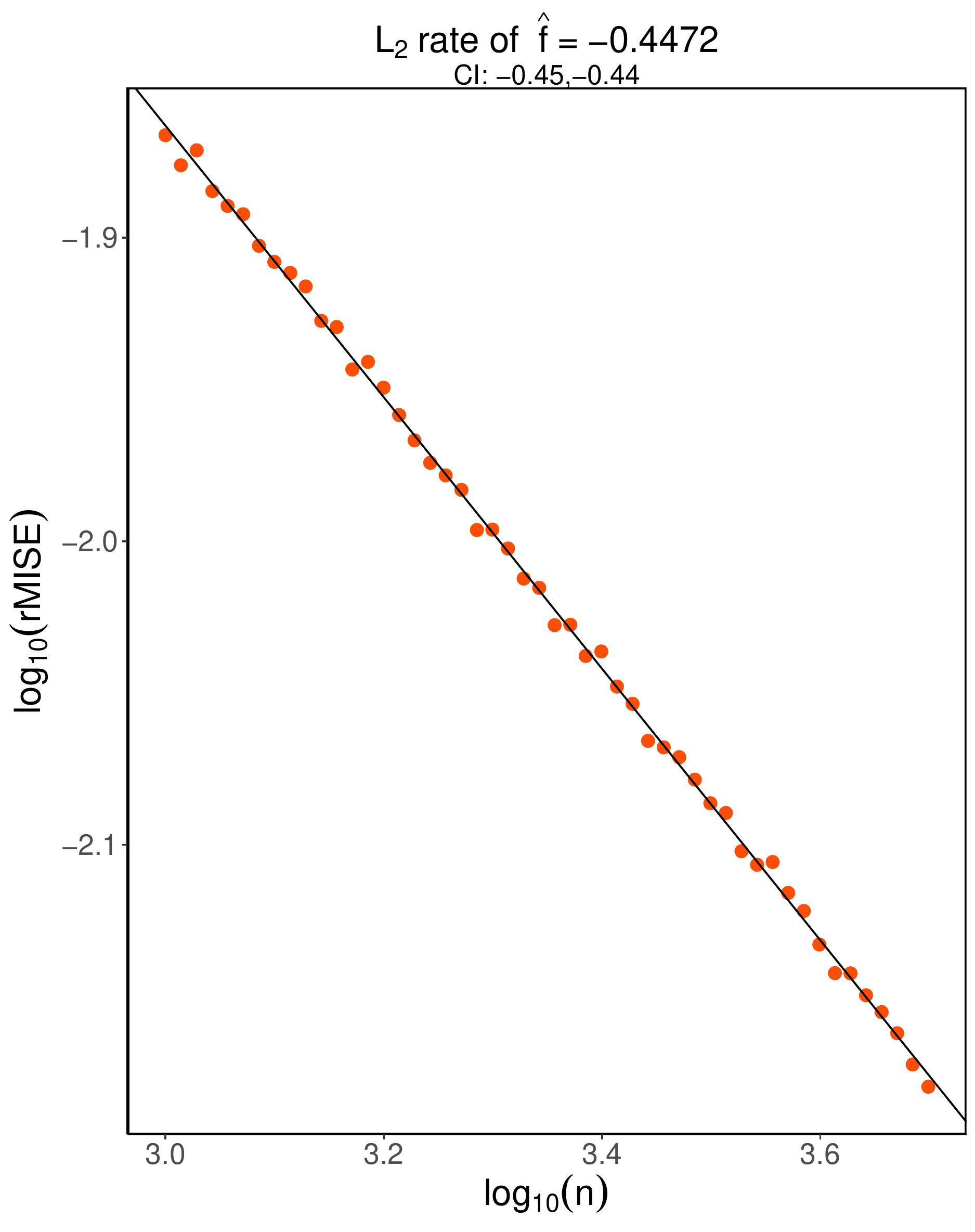} }}%
    \subfloat[\centering \scriptsize $L_2$ rate of convergence for $\hat{f}'$]{{\includegraphics[width=4.6cm]{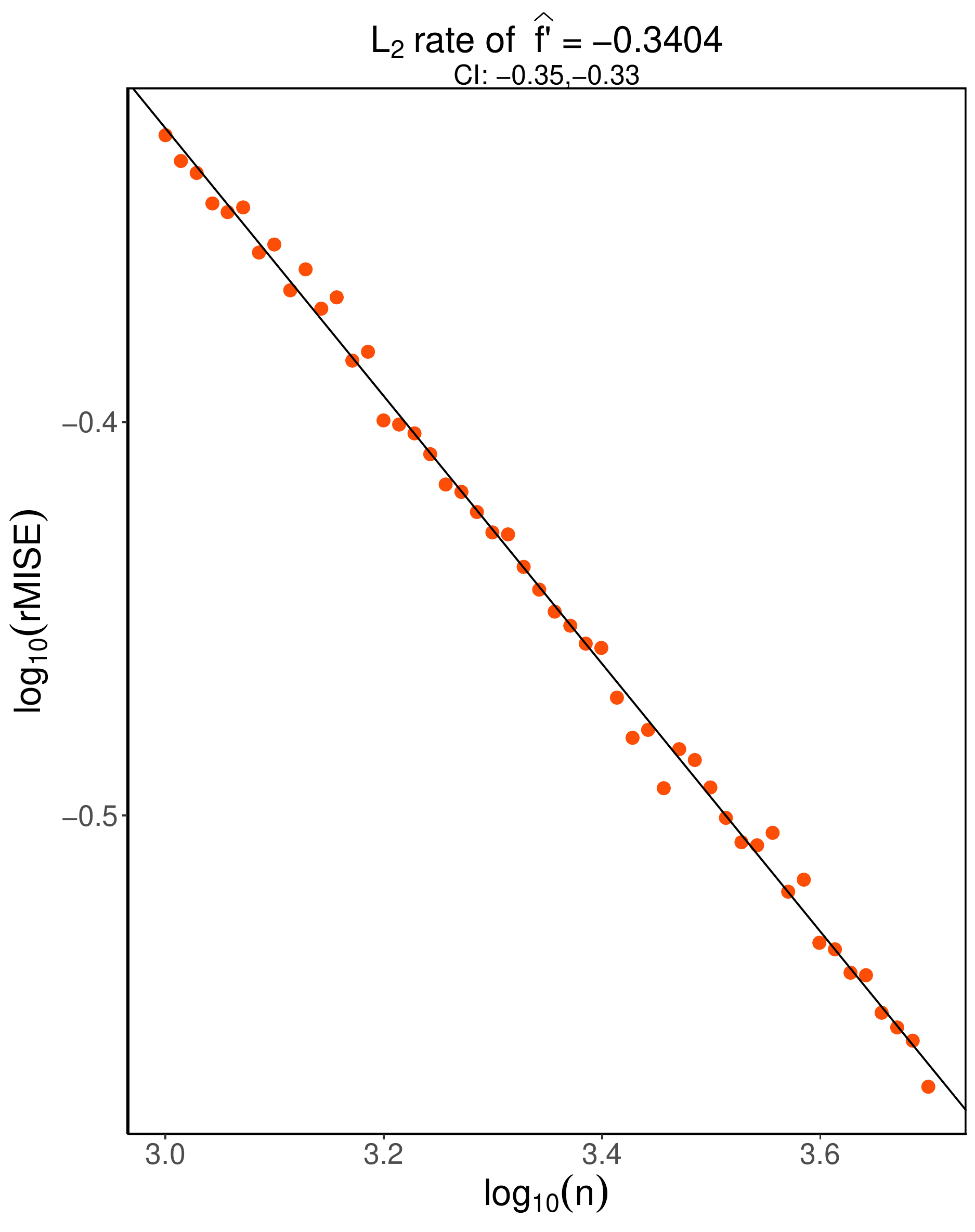} }}%
    \subfloat[\centering \scriptsize $L_2$ rate of convergence for $\hat{f}''$]{{\includegraphics[width=4.6cm]{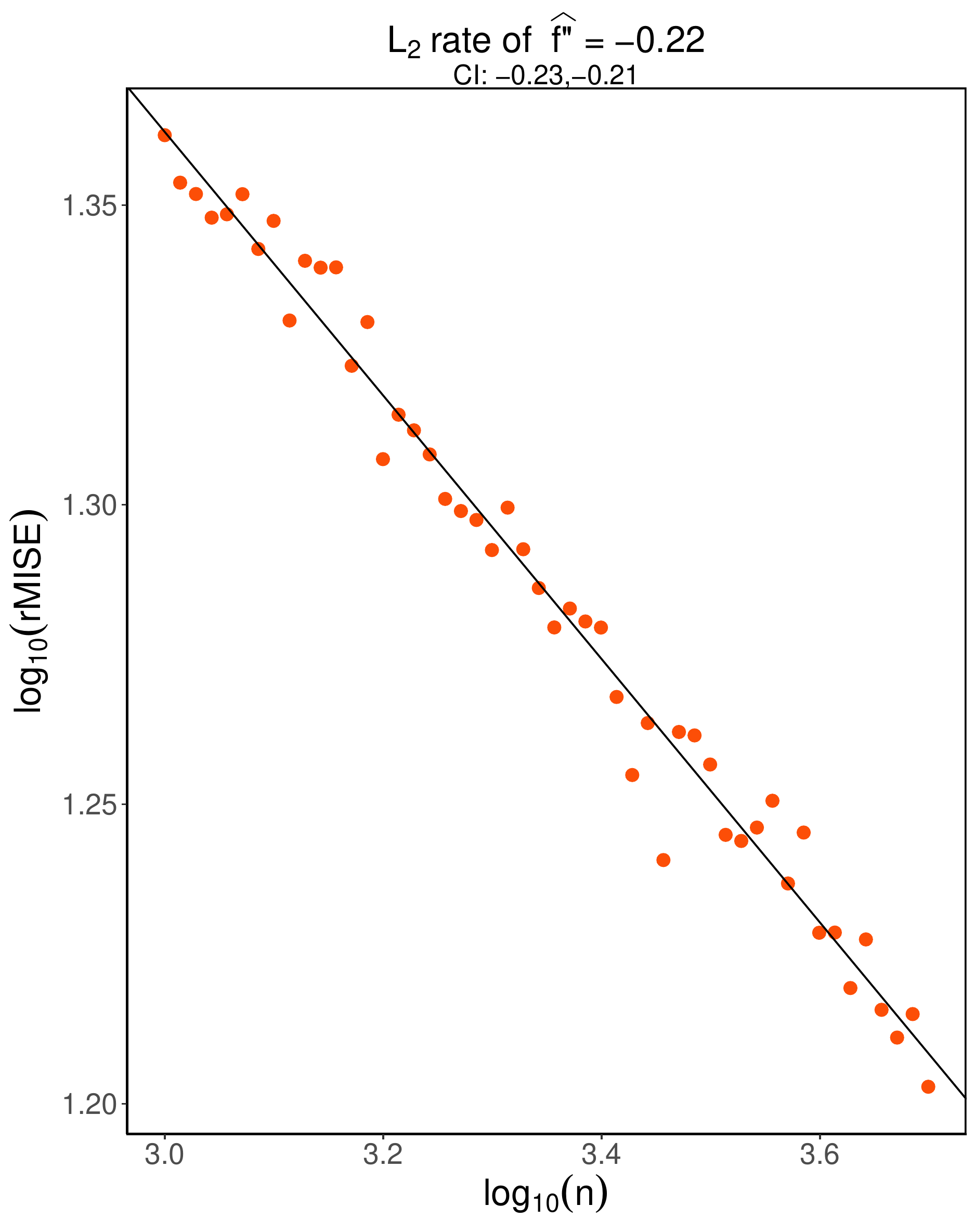} }}%
    \caption{$L_2$ convergence rates for $f$ and its first two derivatives under two scenarios for increasing $K$ with $n$. Figures (a-c) show results for slowly increasing $K$ scenario while Figures (d-f) show results for the fast increasing $K$ scenario. The smoothing parameter $\lambdan$ is chosen by the GCV method.}%
    
    \label{fig:results_gcv}%
\end{figure}

 Figure \ref{fig:results_gcv} illustrates the $L_2$ rate of the naive estimator when the smoothing parameter $\lambdan$ is chosen by the GCV approach. The naive estimator achieves the optimal $L_2$ rates of convergence for the mean regression function and its first two derivatives when GCV is used to choose the smoothing parameter, but it is slightly slower for REML. This deviation from the optimal rate using REML appears to worsen for higher derivatives. Also, we observed that the fast $K$ scenario was overall slightly slower than the slow $K$ scenario for REML. These results agree with known results in the literature for smoothing splines when estimating the mean regression function. For instance, \cite{craven_smoothing_1978} showed that GCV achieves the optimal rate of convergence when used to choose the smoothing parameter in smoothing splines. However, \cite{wahba1985} found that Maximum Likelihood (ML) based methods may be slower than GCV for sufficiently smooth functions.

Table \ref{table:naive} below summarizes the rates of convergence of the naive estimator for estimating derivatives of the mean regression function in \eqref{eqn:barbranter} under the various scenarios of the number of knots $K$ as $n$ increases.
 
 \begin{table}[h!]
    \tiny
     \centering
     \begin{tabular}{ |P{1cm}|P{2cm}|P{2cm}|P{3cm}|P{3cm}| }
         \hline
         \textbf{$\lambdan$ Method} & \textbf{Target} & \textbf{Optimal $L_2$ Rate} & \textbf{Slow $K$} & \textbf{Fast $K$} \\
         \hline
         & $f$   & $-0.44$    &  $-0.45 (-0.45, -0.44)$ & $-0.45 (-0.45, -0.44)$ \\
         GCV & $f'$   & $-0.33$    & $-0.34 (-0.35, -0.33)$ & $-0.34 (-0.35, -0.33)$ \\
         & $f''$   & $-0.22$    & $-0.22 (-0.23, -0.21)$ & $-0.22 (-0.23, -0.21)$ \\
         \hline
         & $f$   & $-0.44$    &  $-0.44 (-0.44, -0.43)$ & $-0.43 (-0.44, -0.43)$ \\
         REML & $f'$   & $-0.33$    & $-0.32 (-0.32, -0.31)$ & $-0.31 (-0.32, -0.31)$ \\
         & $f''$   & $-0.22$    & $-0.19 (-0.19, -0.18)$ & $-0.18 (-0.18, -0.17)$ \\
         \hline
    \end{tabular}
    \caption{Summary of $L_2$ rates of convergence for estimating the mean regression function in \eqref{eqn:barbranter} and its first two derivatives.}
    \label{table:naive}
\end{table}

\subsection{Finite sample performance of naive estimator.}
\label{sec:finitesample}
In this section we compare the naive estimator to an ``oracle'' method that uses knowledge of the true form of the target (mean regression function or its derivatives) to choose the optimal amount of smoothing, which we did with a grid search. While this ``oracle''  is not an estimator, it shows the minimum loss when estimating the function in question with a penalized spline. GCV was used to choose the appropriate smoothing parameter for the various spline-based estimators in what follows.

In Figure \ref{fig:results_median_fits} below, we show the naive estimator That corresponds to the median MSE in the Monte Carlo experiment. To summarize, we see that the naive estimator appears to accurately estimate both the true mean regression function ($f$) and its first derivative ($f'$). However, we observe some lack of fit around the boundaries of the second derivative, ($f''$). 
\begin{figure}[H]
    \centering
    \subfloat[\centering]{{\includegraphics[width=4.5cm]{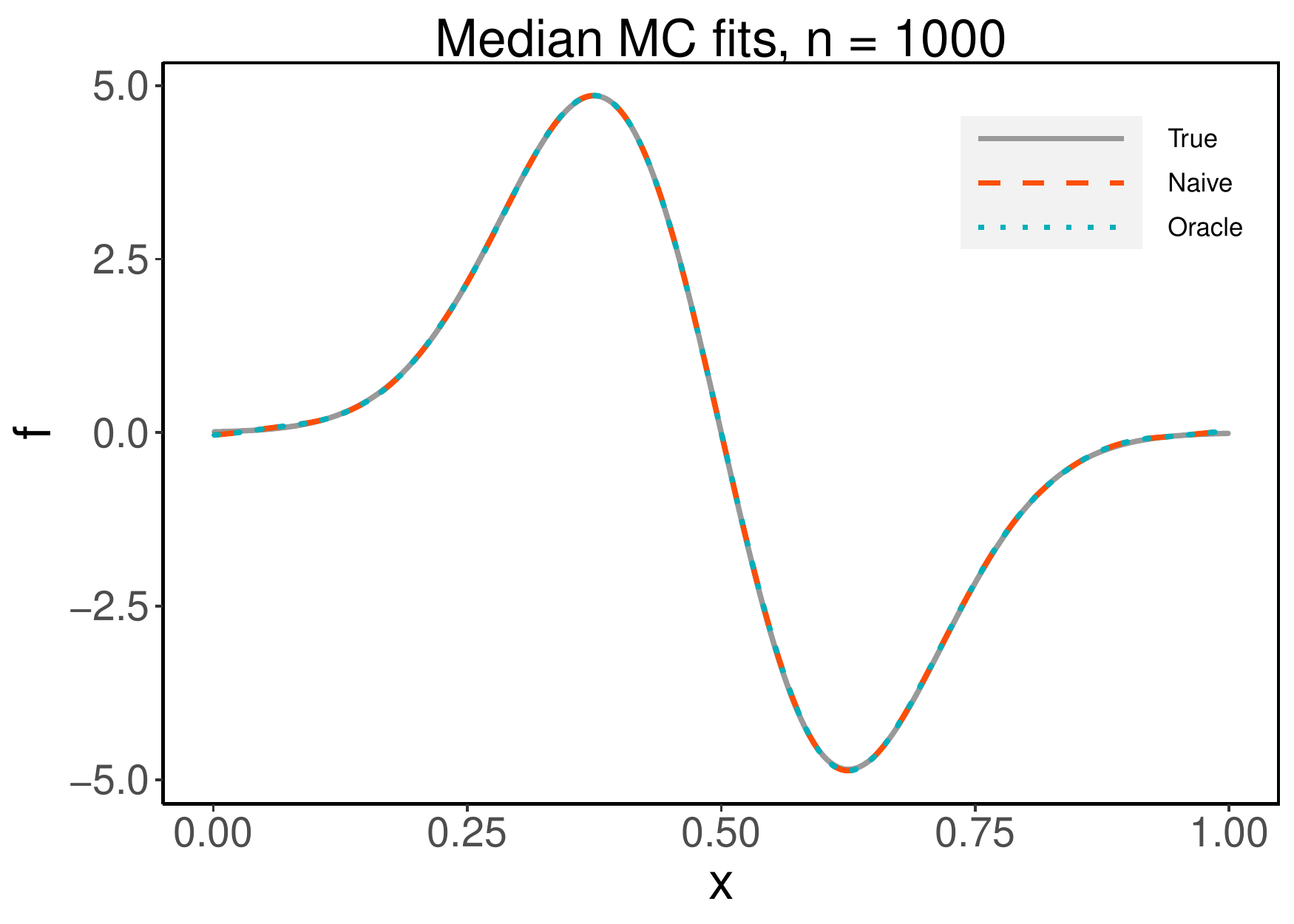} }}%
    \subfloat[\centering]{{\includegraphics[width=4.5cm]{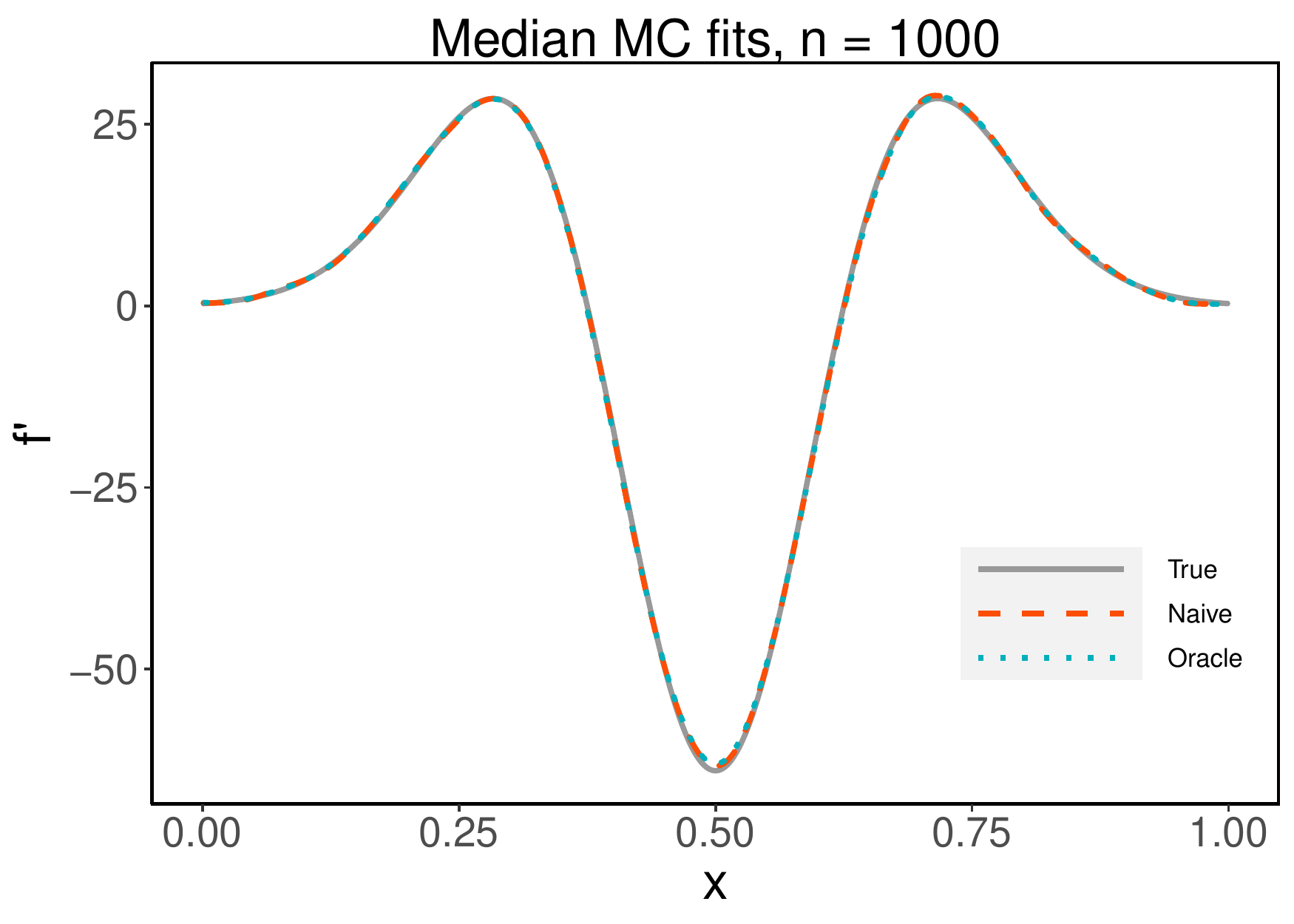} }}%
    \subfloat[\centering]{{\includegraphics[width=4.5cm]{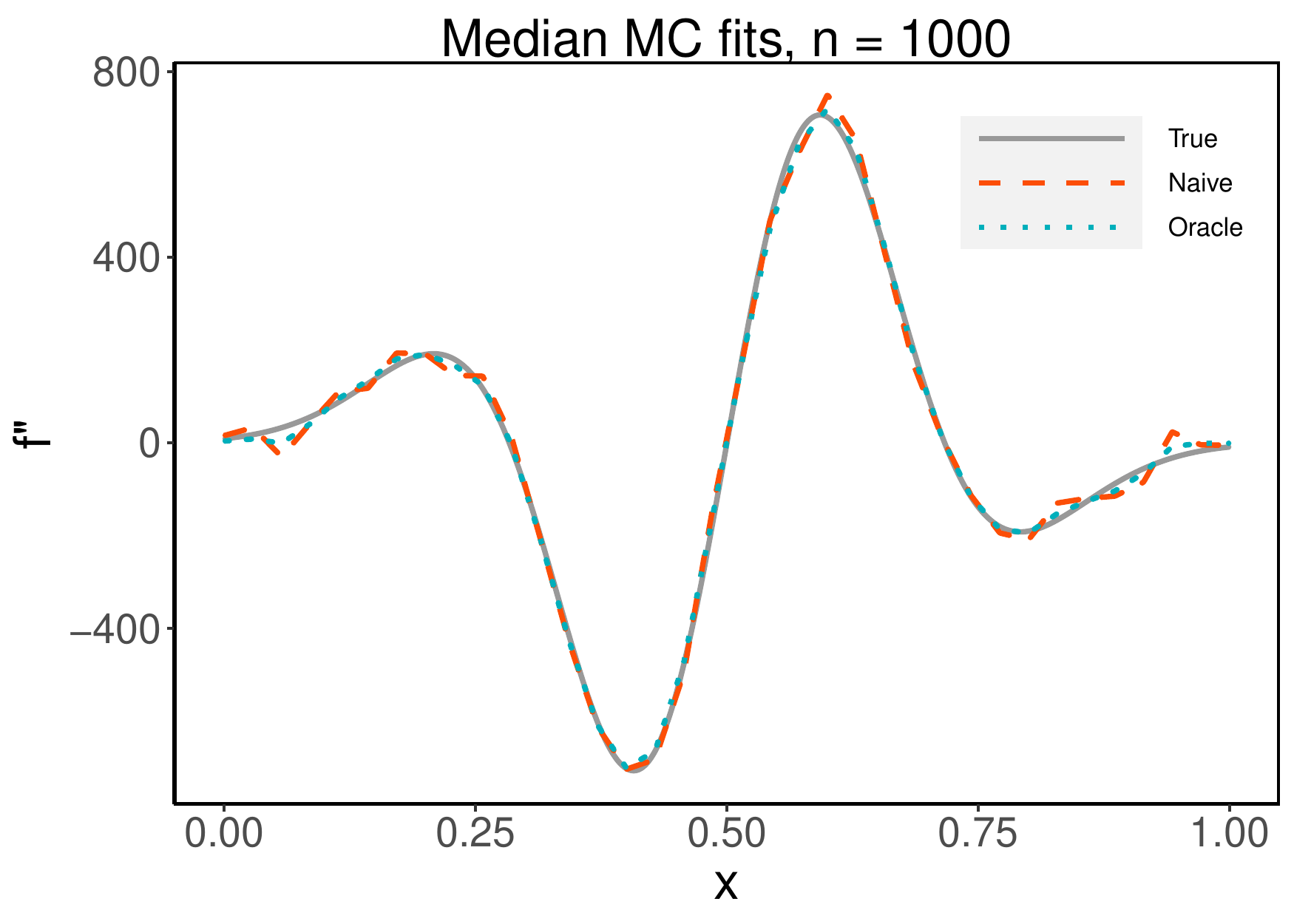} }}%
    \caption{Median Monte-Carlo fits of the mean regression function in \eqref{eqn:barbranter} with its first two derivatives using the naive and oracle estimators.}
    \label{fig:results_median_fits}%
\end{figure}

\begin{figure}[H]
    \centering
    \subfloat[\centering]{{\includegraphics[width=4.2cm]{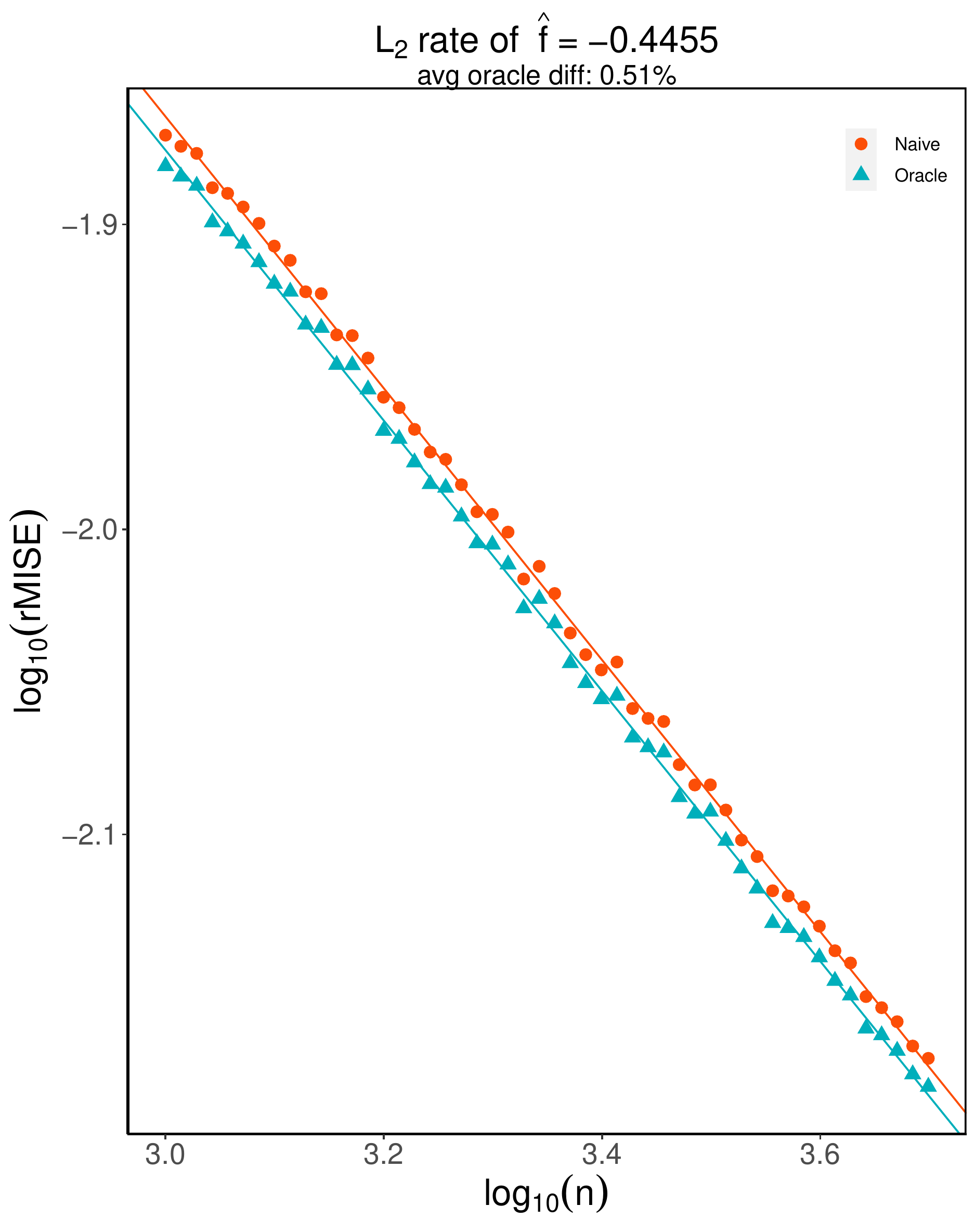} }}%
    \subfloat[\centering]{{\includegraphics[width=4.2cm]{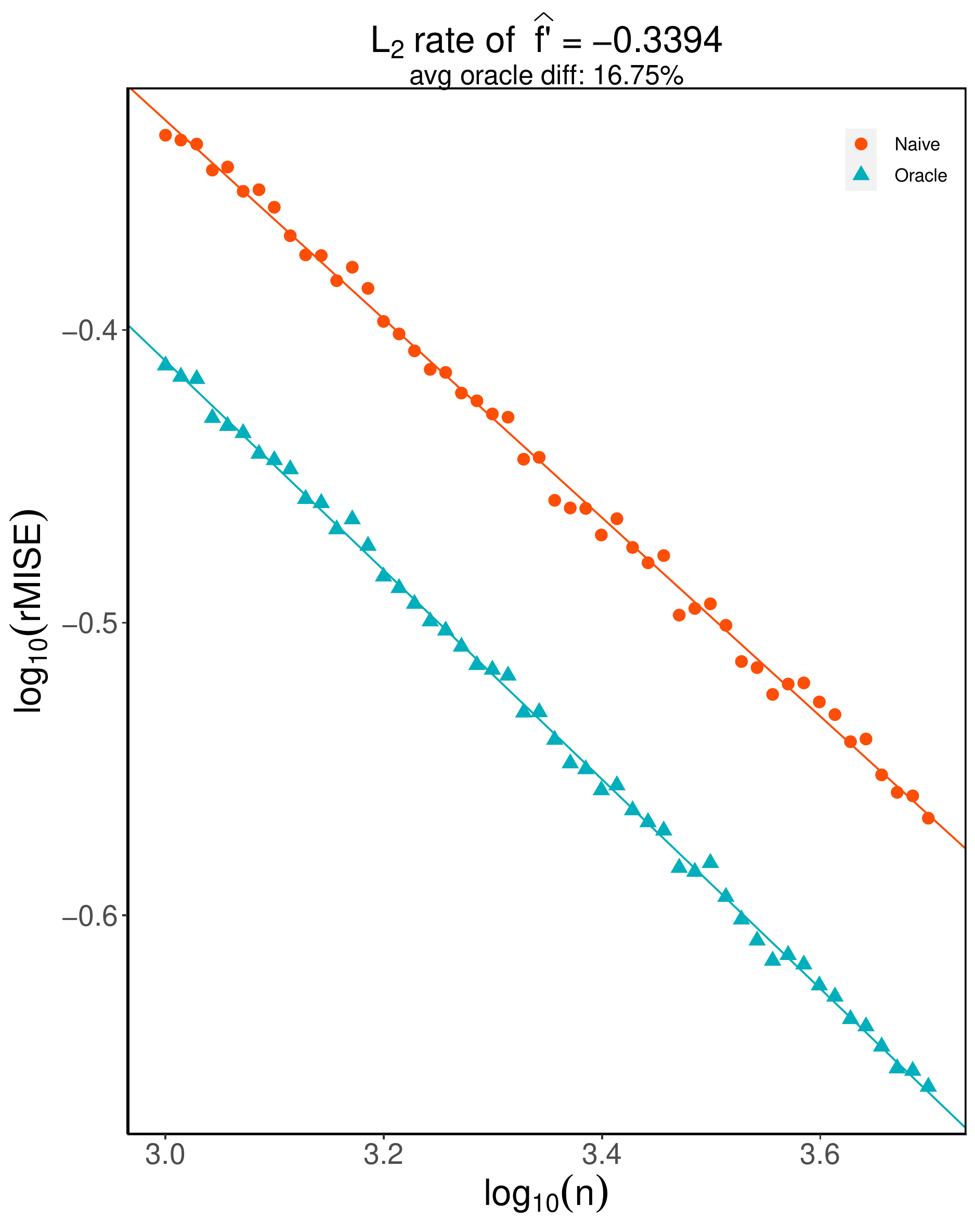} }}%
    \subfloat[\centering]{{\includegraphics[width=4.2cm]{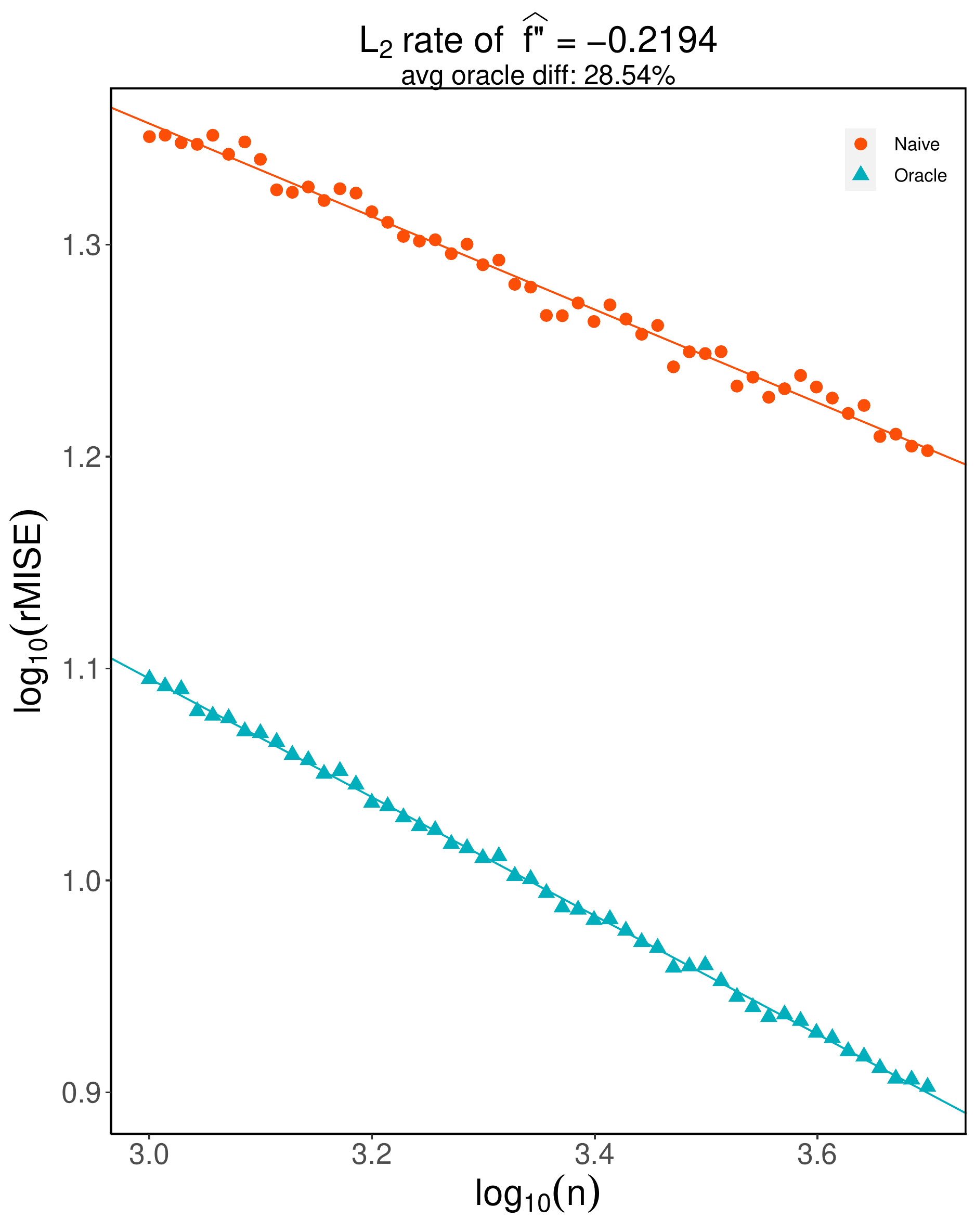} }}%
    
    \centering
    \subfloat[\centering]{{\includegraphics[width=4.2cm]{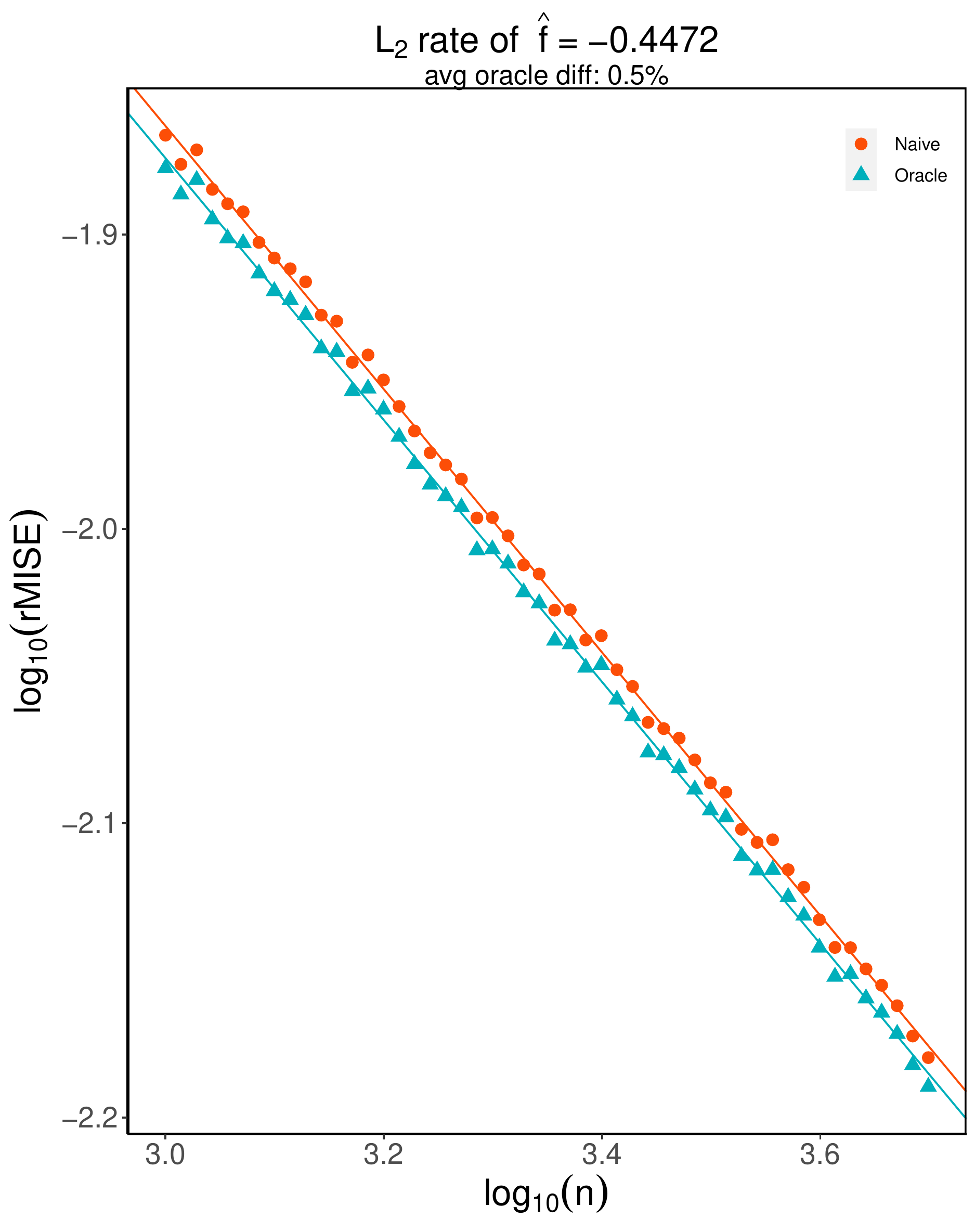} }}%
    \subfloat[\centering]{{\includegraphics[width=4.2cm]{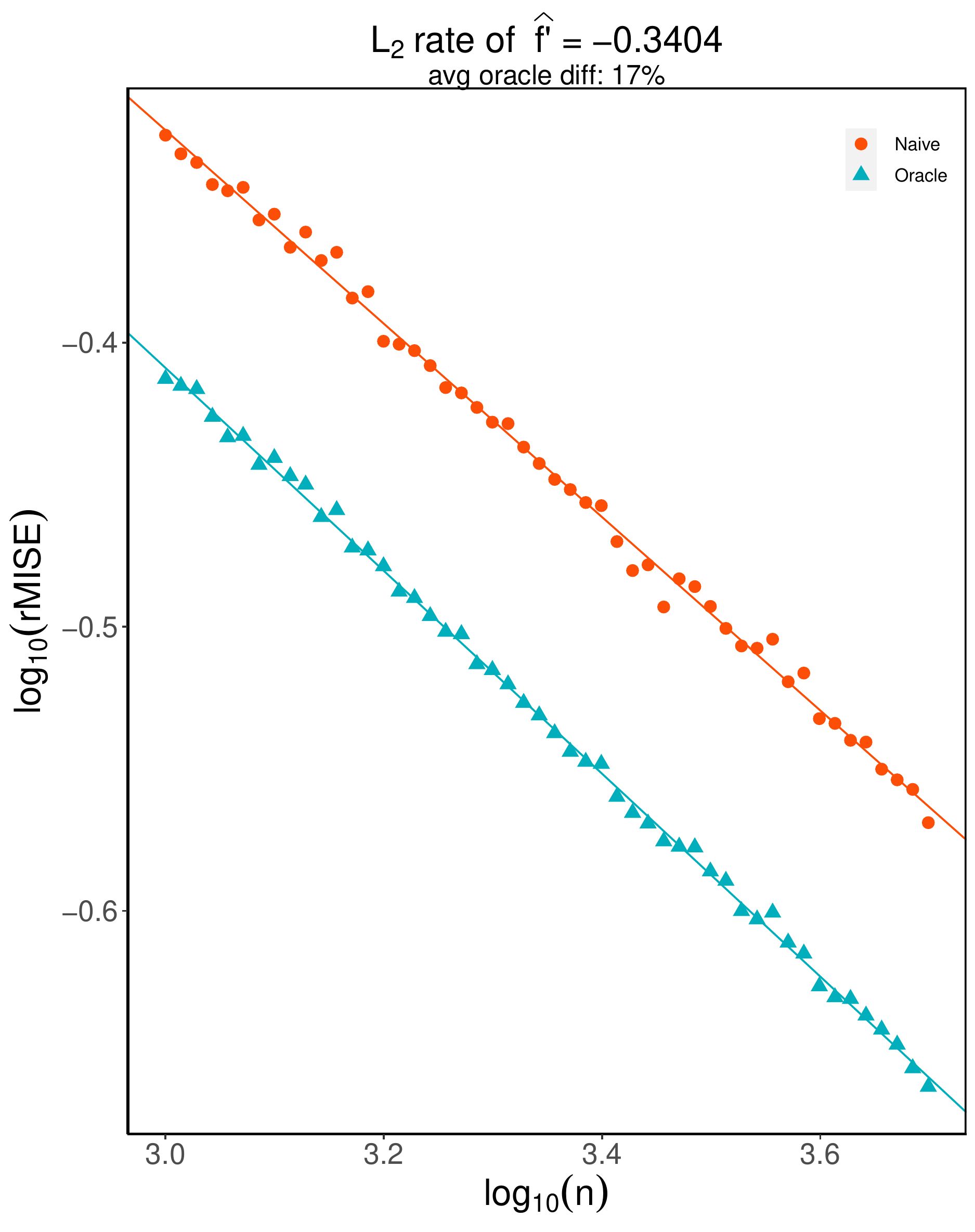} }}%
    \subfloat[\centering]{{\includegraphics[width=4.2cm]{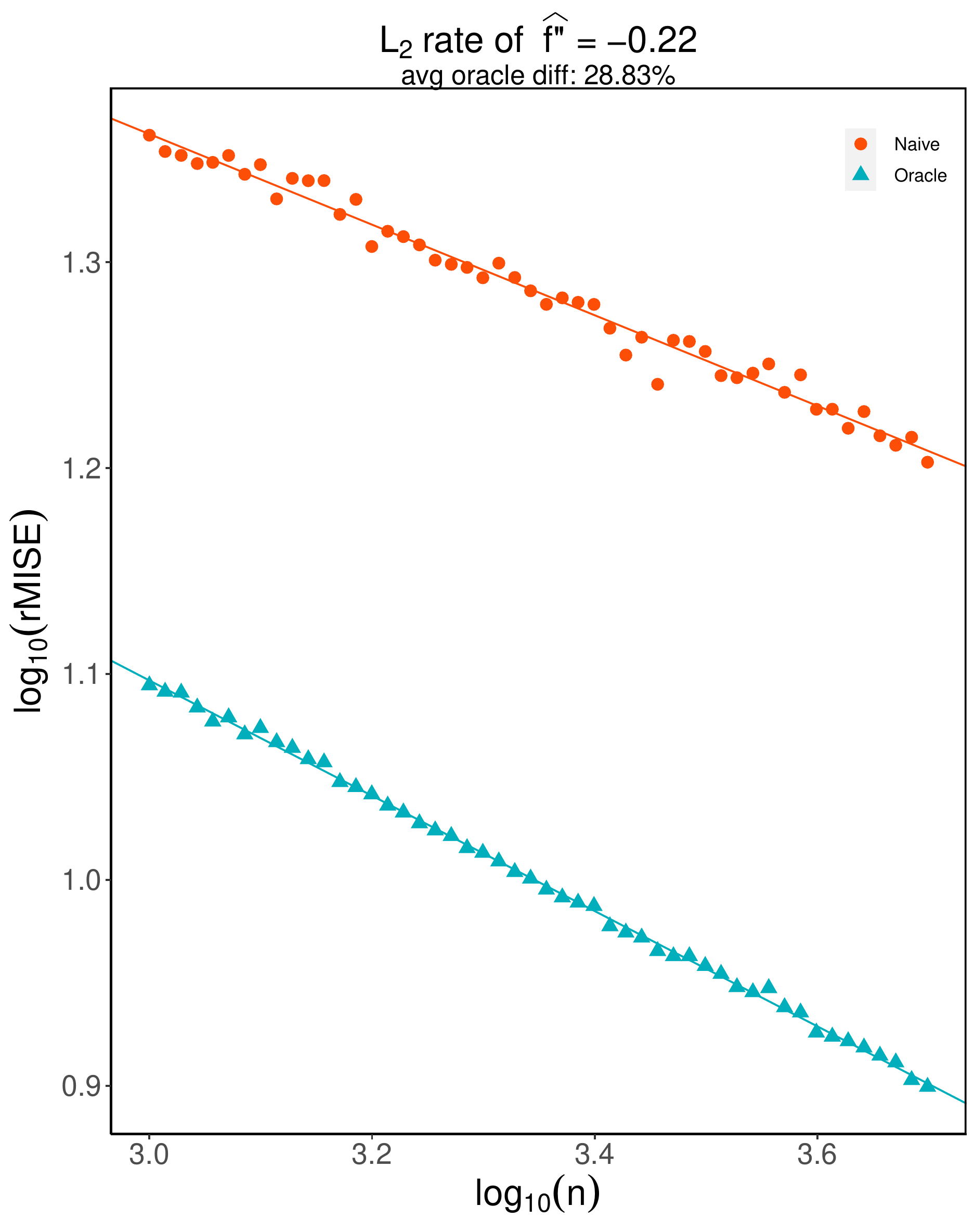} }}%
    \caption{$L_2$ convergence rates for $f$ and its first two derivatives with two scenarios for increasing $K$ with $n$ and how they compare with their corresponding oracle estimators. Figures (a-c) show results for slowly increasing $K$ scenario while Figures (d-f) show results for the fast increasing $K$ scenario. The smoothing parameter $\lambdan$ is chosen by the GCV method.}%
    
    \label{fig:results_oracle}%
\end{figure}

Next, Figure \ref{fig:results_oracle} shows the mean percentage difference between the naive and oracle methods for the mean regression function in \eqref{eqn:barbranter} and its first two derivatives across the two increasing $K$ scenarios. Overall, in comparison to the oracle method, the naive estimator's finite sample performance degrades with increasing derivatives, with an average error difference (logarithmic scale) of about 0.5\% for the mean regression function, 17\% for its first derivative, and 29\% for its second derivative. While the naive penalized spline derivative estimator is shown to converge at the optimal $L_2$ rate of convergence (Theorem \ref{thm:1}), it may also have higher mean squared error in finite samples, especially for higher derivatives. We note that the results summarized in Figure \ref{fig:results_oracle} are similar for the two increasing $K$ scenarios.

\subsection{Comparison with other methods}

In this section, we compare the finite sample MSE of the naive estimator to other derivative estimation methods in the literature. We considered the adaptive penalty penalized spline estimator by \cite{simpkin2013}. We also used the linear combination method of \cite{daitong2016}, but it consistently had higher MSE values and results are not shown. 

We evaluated the methods using three mean regression functions from the literature (\citealp[]{brabanter2013,daitong2016}). As proxies for low, medium, and high noise scenarios, we considered noise levels that were 10 percent, 30 percent, and 60 percent of the range of each function. This was to understand how the methods compare at different levels of noise. The following are the three functions considered:
\begin{equation*}
    f_1(x) = \sin^2(2\pi x) + \log(4/3 + x) \quad \textrm{for} \quad x\in[-1, 1],
\end{equation*}
\begin{equation*}
    f_2(x) = 32 e^{-8(1-2x)^2}(1-2x) \quad \textrm{for} \quad x\in[0, 1],
\end{equation*}
and the doppler function
\begin{equation*}
    f_3(x) = \sqrt{x(1-x)} \sin\left(\frac{2.1\pi}{x+0.05}\right) \quad \textrm{for} \quad x\in[0.25, 1].
\end{equation*}

Figure \ref{fig:results_comparisons} below shows the results for estimating the first (panel a) and second (panel b) derivatives of the three mean regression functions across the three noise levels. These results indicate that the adaptive penalty methods and the naive method often perform similarly, depending on the form of the function, the noise level, and the order of the derivative. We also note that the adaptive penalty method sometimes performs better than the oracle method. This is possible since the oracle method only finds the best P-splines estimate with the form of the penalty held constant. 

\begin{figure}[H]
    \centering
    \subfloat[\centering]{{\includegraphics[width=7cm]{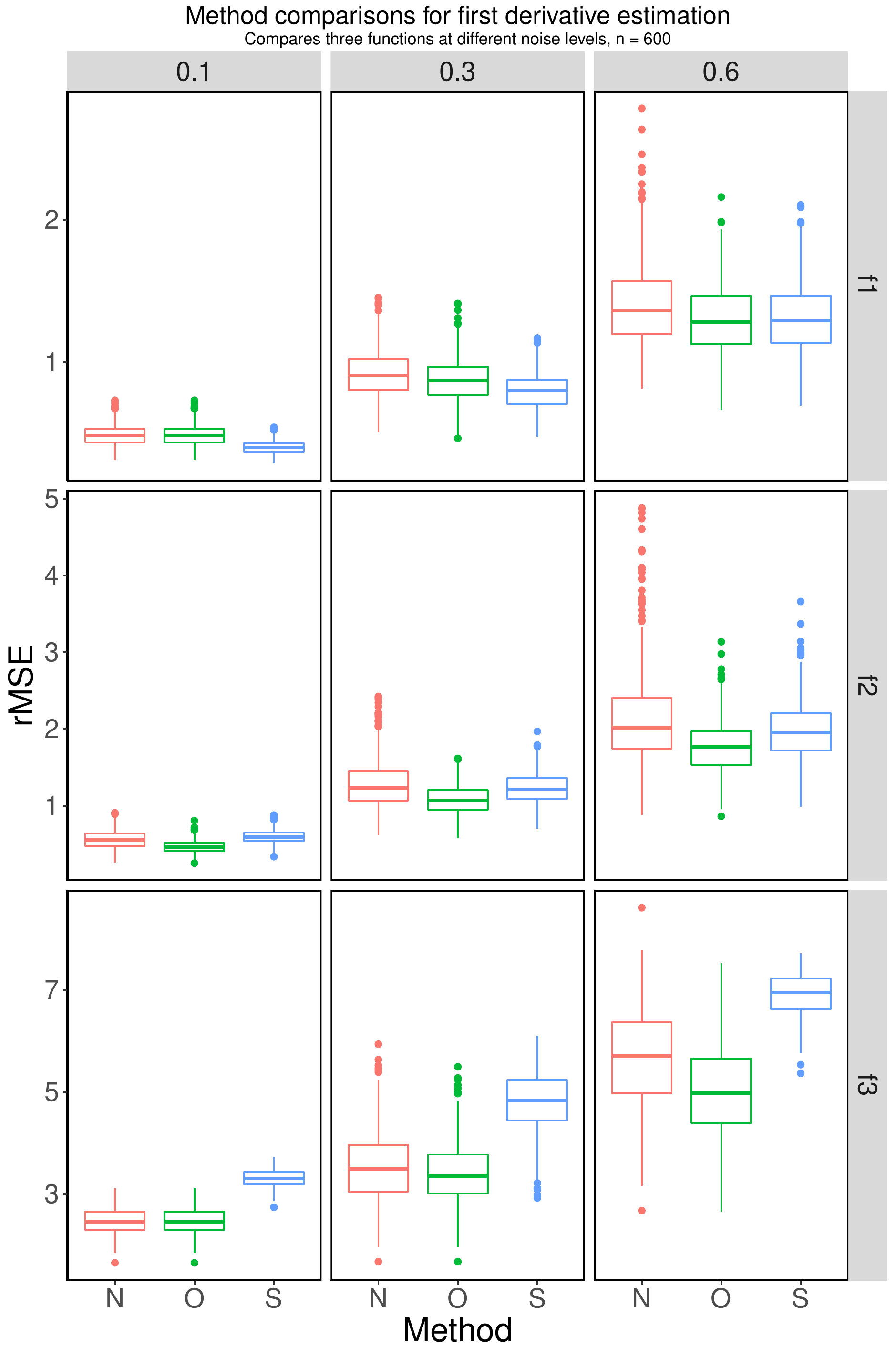} }}%
    \subfloat[\centering]{{\includegraphics[width=7cm]{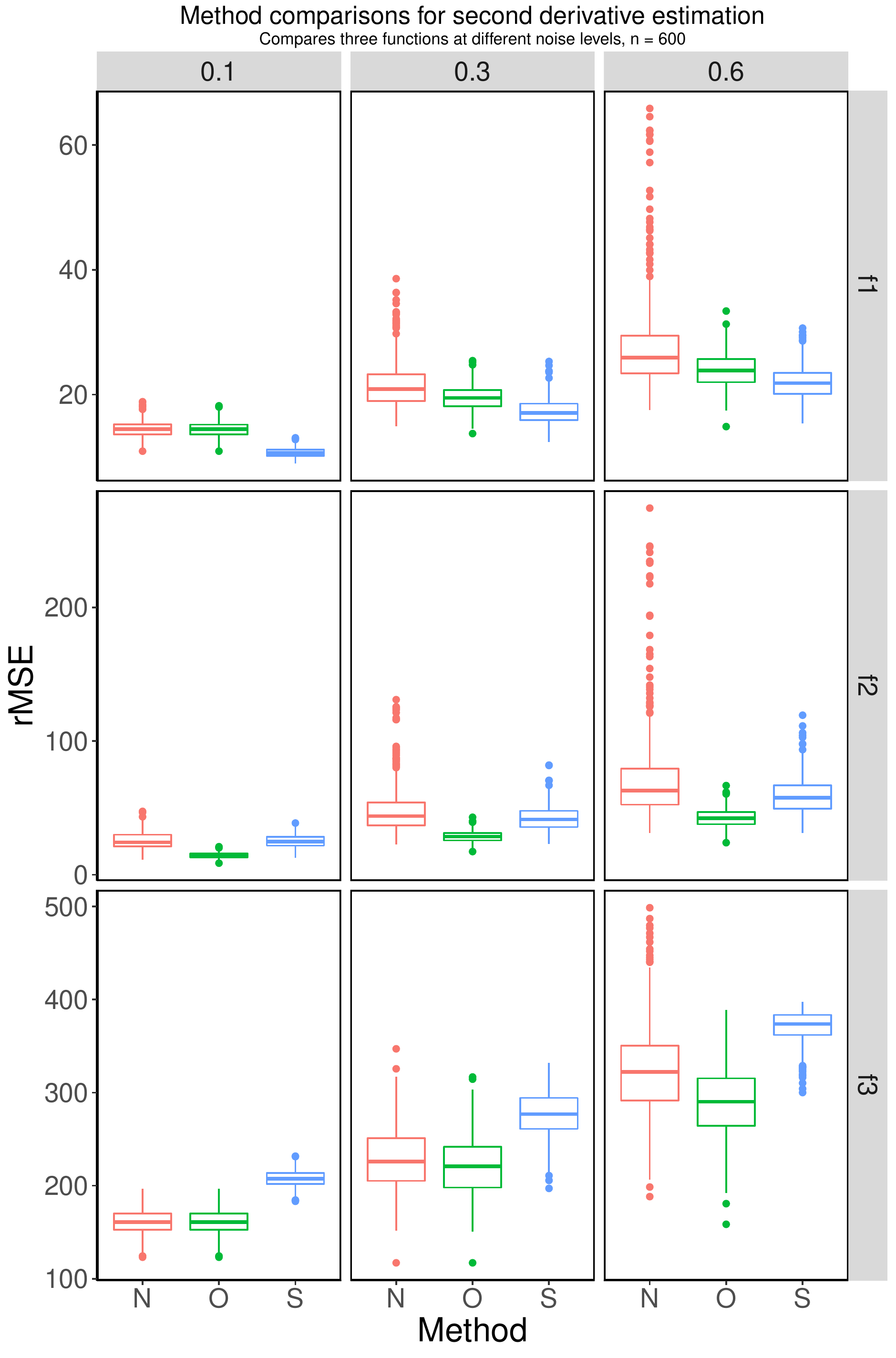} }}%
    \caption{Comparing derivative estimation methods in reference to the oracle estimator across different functions and noise levels. Panel (a) shows results for estimating first derivatives of the mean regression functions $f_1$, $f_2$ and $f_3$ while Panel (b) shows results for estimating second derivatives.}%
    \label{fig:results_comparisons}%
\end{figure}

%% file: conclusion.tex
In this paper, we have shown that the naive penalized spline estimator of the $r^{th}$ derivative of the mean regression function achieves the optimal $L_2$ rate of convergence (\citealp[]{stone1982}) under standard assumptions on knot placement and the penalty matrix. This builds on the work by \cite{xiao2019} which derived the $L_2$ rate of convergence for estimating the mean regression function. As stated in Remark \ref{rmk:1} and noted by others (\citealp[]{Claeskens2009, xiao2019}), the rate at which the number of knots, $K$, increases with $n$ gives rise to two scenarios: the fast $K$ scenario, which is similar to smoothing spline asymptotics, and the slow $K$ scenario, which is similar to regression spline asymptotics.

Using simulations, we investigated how two prevalent methods for choosing the smoothing parameter (GCV and REML) affect the $L_2$ convergence of the naive estimator. We found that, for both slow and fast K scenarios, the naive estimator  achieves the optimal $L_2$ rate of convergence when GCV is used. For REML, the estimator did not quite achieve the optimal rate. 

To access the finite sample performance of the naive estimator, we compared the MSEs of the estimator with an ``oracle'' method that uses information about the true function to be estimated to choose the P-spline's smoothing parameter. We found that, in finite samples, the naive estimator may have noticeably larger mean squared errors, especially for higher derivatives, but the estimates can still be quite visually similar. We found that the adaptive penalty penalized spline estimator by \cite{simpkin2013} performed similarly to the naive estimator.

%% file: appendix.tex
\subsection{Proof of Theorem}
The proof proceeds in two steps. We first derive the $L_2$ rate of convergence for the bias of the naive estimator and then we derive that of the variance. The approach of the proof closely follows the proof for the $L_2$ rate of convergence of the mean regression function itself found in \cite{xiao2019}.
We start by defining some terms to simplify the notation. \\
Let $G_{n, q} = \frac{\bt\bb}{n}$ and $H_n = G_{n, q} + \lambdan \Pm$.
To ease exposition, we follow \cite{zhou2000yhh} and write $\fhp$ as
$$\fhp = \bqq D^{(r)}(\gn+\lambdan\Pm)^{-1}\bt\bm{y}/n$$
where $\bqq \in \mathbb{R}^{K+q-r}$ is a vector of B-Spline basis functions of order $q-r$ and $\dr$ is defined as $\dr = M^T_r\times M^T_{r-1}\times \dots \times M^T_1$ with
$$
M_l = (q-1)\begin{bmatrix} 
    \frac{-1}{t_1-t_{1-q+l}} & 0 & 0 & \dots & 0 \\
    \frac{1}{t_1-t_{1-q+l}} & \frac{-1}{t_2-t_{2-q+l}} & 0 & \dots & 0 \\
    0 & \frac{1}{t_2-t_{2-q+l}} & \frac{-1}{t_3-t_{3-q+l}} & \dots & 0 \\
    \vdots & \vdots & \vdots & \vdots & \vdots\\
    0 & 0 & 0 & \dots & \frac{1}{t_{K+q-l}-t_K} \\
    \end{bmatrix}
$$
for $1\le l\le r\le q-2$.

Let $\bpq=\bqq D^{(r)}$, implying
$$ \hat{f}^{(r)}(x) = \bpq\left(G_{n, q}+\lambdan \Pm\right)^{-1}\bt\bm{y}/n$$
We use the identity $(A+B)^{-1} = A^{-1} - A^{-1}B(A+B)^{-1}$ to expand the inverse term in the estimator. This later allows us to split the bias term into the part due to approximating $\fp$ with a spline (approximation bias) and the other part due to penalization (shrinkage bias).
\begin{eqnarray*}
(G_{n, q}+\lambdan \Pm)^{-1} &=& G_{n,q}^{-1}-G_{n,q}^{-1}(\lambdan \Pm)H_n^{-1} \\
&=& G_{n,q}^{-1}-H_n^{-1}(\lambdan \Pm)G_{n,q}^{-1}
\end{eqnarray*}
where the last equality is by symmetry.\\
Substituting into $\fph$, we have:
\begin{eqnarray*}
    \implies \hat{f}^{(r)}(x) &=& \bpq\left(\Gni - \Hni(\lambdan \Pm)\Gni\right)\Bt\y/n\\
    &=& \bpq\Gni\Bt\y/n - \bpq\Hni\lambdaP\Gni\Bt\y/n
\end{eqnarray*}
We now focus on the bias of the naive estimator, $\efph - \fp$.

From Lemma \ref{lm1}, 
$\exists s_f\in \Sqt$, the space of spline functions of order $q$ defined on knots $\underline{\mathbf{t}}$ such that
$||f^{(r)} - s_f^{(r)}|| = O(h^{q-r}) + o(h^{p-r}).$
The bias of $\fhp$ can be written as:
\begin{equation}
\label{eqn:bias}
    \efph - \fp = \left[E\left(\fph\right) - \sfp\right] + \left[\sfp - \fp\right]
\end{equation}
Equation \eqref{eqn:bias} above allows us to separately evaluate the approximation bias and shrinkage bias for estimating $\fp$.
Notice that Lemma \ref{lm1} provides information on the rate of convergence on the second term in \eqref{eqn:bias}, we will next focus expressing the first term in a form that isolates the effect of penalization on the bias.  Substituting the previously derived expression for $\fph$ into the first term, we have
\begin{eqnarray*}
 E\fph - \sfp &=& \bpq\Gni\Bt\f/n - \sfp \\
 && - \bpq\Hni\lambdaP\Gni\Bt\f/n  \\
 &=& \bpq\gammab - \sfp - \bpq\hni\lambdap\gammab
\end{eqnarray*}
where $\gammab = \Gni\Bt\f/n$ and $\f = E\left[\y\right]$.

But $\sfp = \bpq\Gni\Bt\sfb/n$, where $\sfb = \left\{s_f(x_1), s_f(x_2), \dots, s_f(x_n)\right\}$
\begin{eqnarray}
\label{eqn:pen_bias}
\implies E\fph - \sfp &=& \bpq\gammab - \bpq\Hni\lambdaP\gammab \nonumber \\
&& - \bpq\Gni\Bt\sfb/n \nonumber\\
&=& \bpq\Gni\Bt\f/n - \bpq\Gni\Bt\sfb/n \nonumber\\ 
&& - \bpq\Hni\lambdaP\gammab \nonumber\\
&=& \bpq\Gni\Bt(\f-\sfb)/n - \bpq\Hni\lambdaP\gammab \nonumber\\
&=& \bpq\Gni\alphab - \bpq\hni\lambdaP\gammab
\end{eqnarray}
where $\alphab = \Bt(\f-\sfb)/n$

Let $Q(x)$ be a distribution of $x$ on $\left[0, 1\right]$ with positive continuous density $q(x)$. Then substituting \eqref{eqn:pen_bias} into \eqref{eqn:bias} and using the triangle inequality, we can evaluate the squared bias of $\fhp$ as:
\begin{eqnarray}
\label{eqn:bias_l2}
\int_0^1 \left(\efph - \fp\right)^2 q(x)dx &\le& \int_0^1\left(\sfp - \fp\right)^2q(x)dx \nonumber\\
 && + \alphabt\gni\gp\gni\alphab\\
 &&  + \gammabt\lambdap\hni\gp\hni\lambdap\gammab  \nonumber 
\end{eqnarray}
where $\gp=\int_0^1 \bpqt\bpq q(x)dx$. The first and second terms in \eqref{eqn:bias_l2} represent the part of the bias due to using spline functions to estimate $\fp$, and the last term represents the part of the bias due to penalization.

Observe that, by Lemma \ref{lm1},
\begin{eqnarray*}
\int_0^1\left(\sfp - \fp\right)^2q(x)dx &\le& q_{\max} \int_0^1 \left(\sfp - \fp\right)^2dx \\
&=& O\left(h^{2(q-r)}\right) +  o\left(h^{2(p-r)}\right)
\end{eqnarray*}
where, $q_{\max} = \displaystyle \max_{0\le x\le 1}q(x) < \infty$.

For the second term in \eqref{eqn:bias_l2}, we use the result $||\gp||_\infty = O(h^{-2r})$, from Lemma \ref{lm2}. We also use $||\gni||_\infty=O(h^{-1})$ from Lemma \ref{lm4} and Lemma 6.10 of \cite{agarwal1980} that $||\alphab||_{\max}=o(h^{p+1})$.

Let $\gph$ be square and symmetric matrix such that $\gp=\gph\gph$.

We write 
\begin{eqnarray*}
\alphabt\gni\gp\gni\alphab &=& \left(\gph\gni\alphab\right)^T\left(\gph\gni\alphab\right)\\
&=& ||\gph\gni\alphab||_2^2 \\
&\le& ||\alphab||_2^2||\gph\gni||_2^2\\
&\le& ||\alphab||_2^2||\gph||_2^2||\gni||_2^2\\
&=& o(h^{2p+2})O(h^{-2r})O(h^{-2}) \\
&=& o(h^{2p-2r})
\end{eqnarray*}
Next, we focus on the part of the bias due to penalization as given by the third term in \eqref{eqn:bias_l2}. First, note that from \cite{deboor1978} and Lemma 5.2 of \cite{zhou2000yhh}, $\dr$ in $\bpq=\bqq D^{(r)}$, is such that
$$||D^{(r)}||_\infty = O(h^{-r})$$ 
This can be easily seen by inspecting the elements of $\dr$.
$$\therefore \bpqt\bpq=D^{(r)T}\bqqt \bqq\dr = O(h^{-2r})\bqqt\bqq$$

Thus, we can write
\begin{eqnarray*}
\gp&=&\int_0^1 \bpqt\bpq q(x)dx \\
&=& O(h^{-2r})\int_0^1 \bqqt\bqq q(x)dx \\
&=& O(h^{-2r}) G_{q-r}
\end{eqnarray*}
where $G_{q-r} = \int_0^1 \bqqt\bqq q(x)dx$.

Also, by the WLLN, $G_{n, q-r}=G_{q-r}+o(1)$.\\
Therefore:
\begin{eqnarray*}
\psib &=& O(h^{-2r})\gammabt\lambdap\hni \\ 
&& \times\ G_{q-r}\hni\lambdap\gammab\\
&=& O(h^{-2r})\gammabt\lambdap\hni \\ 
&& \times\ \gnqq\hni\lambdap\gammab
\end{eqnarray*}
where $\gnqq = B_{q-r}^TB_{q-r} / n$, the version of $\gn$ based on B-splines of degree $q-r$.
Note that the decay of the eigenvalues of $G_{n, q}$ does not depend on $q$ (see Lemma \ref{lm3}). Therefore, we will use
$G_{n, q}$ instead of $G_{n, q-r}$ in the derivations that follow for asymptotic order. This simplifies the calculations since $\hni$ depends on $\gn$.

Using 
\begin{eqnarray*}
\hni&=&\left[\gn+\lambdap\right]^{-1} \\
&=& \left[\gnh\left(\gnh+\lambdan\gnmh \Pm\right)\right]^{-1}\\
&=& \left(\gnh+\lambdan\gnmh \Pm\right)^{-1}\gnmh
\end{eqnarray*}
we can write
\begin{eqnarray}
\label{eqn:3inner}
\lambdap\hni\gn\hni\lambdap &=& \lambdap\left(\gnh+\lambdan\gnmh \Pm\right)^{-1}\gnmh\gn \nonumber \\ 
&& \times\ \left(\gnh+\lambdan\gnmh \Pm\right)^{-1}\gnmh\lambdap \nonumber
\end{eqnarray}

Let $\tilde{P} = \tildp$ $\implies \tilde{P}\gnh = \gnmh\lambdap$

Substituting into \eqref{eqn:3inner}, we have
\begin{eqnarray*}
\lambdap\hni\gn\hni\lambdap &=& \lambdap\left(\gnh+\tilde{P}\gnh\right)^{-1}\gnh \\
&& \times\ \left(\gnh+\tilde{P}\gnh\right)^{-1}\tilde{P}\gnh\\
&=& \lambdap\gnmh\left(I+\tilde{P}\right)^{-1}\gnh\gnmh \\ && \times\ \left(I+\tilde{P}\right)^{-1}\tilde{P}\gnh\\
&=& \gnh\tilde{P}(I+\tilde{P})^{-2}\tilde{P}\gnh
\end{eqnarray*}
where in the second equality, we've used the fact that $\gnh+\tilde{P}\gnh = (I+\tilde{P})\gnh$
and that $\lambdap\gnmh=\gnh\tilde{P}$ in the last equality.

Using the above, we can then write:
\begin{eqnarray*}
    \gammabt\lambdap\hni\gn\hni\lambdap\gammab &=& \gammabt\gnh\tilde{P}\left(I+\tilde{P}\right)^{-2}\tilde{P}\gnh\gammab
\end{eqnarray*}

From the fact that $\tilde{P}\left(I+\tilde{P}\right)^{-2}\tilde{P}\le ||\tilde{P}||_2\tilde{P}$,
\begin{eqnarray*}
\psib &=& O(h^{-2r})||\tilde{P}||_2\gammabt\gnh\tilde{P}\gnh\gammab \\
&=& O(h^{-2r})||\tilde{P}||_2\gammabt\lambdap\gammab\\
&=& O(h^{-2r})||\gnmh\lambdap\gnmh||_2\gammabt\lambdap\gammab \\
&=& O(h^{-2r})||\gni||_2||\lambdap||_2 \gammabt\lambdap\gammab \\
\end{eqnarray*}
where we have used $\gnh\tilde{P}\gnh=\lambdap$ in the second equality and substituted $\tilde{P}$ in the third.

By Assumption 5, $||\Pm||_2=O(h^{1-2m})$ and from Lemma 6, $\gammabt P_m\gammab = O(1)$.\\
Therefore:
\begin{eqnarray*}
    \psib &=& O(h^{-2r})O(h^{-1})O(\lambdan h^{1-2m})O(\lambdan) \\
    &=& O(\lambdan^2h^{-2m-2r}).
\end{eqnarray*}
Also, from $\tilde{P}(I+\tilde{P})^{-2}\tilde{P}\le \tilde{P}$, we have
\begin{eqnarray*}
\psib &=& O(h^{-2r})\gammabt\gnh\tilde{P}(I+\tilde{P})^{-2}\tilde{P}\gnh\gammab \\
&=& O(h^{-2r})\gammabt\gnh\tilde{P}\gnh\gammab\\
&=& O(h^{-2r})\gammabt\lambdap\gammab \\
&=& O(\lambdan h^{-2r})
\end{eqnarray*}

$\therefore \psib = O\left\{\min\left(\lambdan^2h^{-2m-2r}, \lambdan h^{-2r}\right)\right\}$\\
This concludes the proof for bias in \eqref{eqn:bias_l2}.

Next, we look at the variance part:
\begin{eqnarray*}
Var(\fph) &=& Var\left(\bpq\hni\bt\y/n\right)\\
&=& \bpq\hni\bt Var(\y/n) \bb\hni\bp \\
&=& \frac{\sigma^2}{n} tr\left\{\bpq\hni\bt \bb/n\hni\bp\right\} \\
&=& \frac{\sigma^2}{n} tr\left\{\hni\gn\hni\bp\bpq\right\}
\end{eqnarray*}
\begin{eqnarray*}
\int_0^1 Var(\fph)q(x)dx &=& \frac{\sigma^2}{n}tr\left\{\hni\gn\hni\gp\right\} \\
&=& O(h^{-2r})\frac{\sigma^2}{n}tr\left\{\hni\gn\hni\gn\right\} 
\end{eqnarray*}

From 
\begin{eqnarray*}
\hni &=& \left(\gn+\lambdap\right)^{-1} \\
&=& \left[\gn\left(I+\gni\lambdap\right)\right]^{-1}\\
&=& \left[I+\gni\lambdap\right]^{-1}\gni
\end{eqnarray*}
$$\implies \hni\gn = \left[I+\gni\lambdap\right]^{-1}.$$
Note that $\gni\lambdap = \gnmh\gnmh\lambdap$ and by the rotation property of the trace,
\begin{eqnarray*}
tr\left[\gni\lambdap\right] &=& tr\left[\gnmh\gnmh\lambdap\right] \\
&=& tr\left[\gnmh\lambdap\gnmh\right]\\
&=& tr\left[\tilde{P}\right]
\end{eqnarray*}
\begin{eqnarray*}
\therefore \int_0^1 Var(\fhp)q(x)dx &=& O(h^{-2r})\frac{\sigma^2}{n}tr\left[(I+\tilde{P})^{-2}\right]\\
&=& O(h^{-2r})\frac{\sigma^2}{n}||(I+\tilde{P})^{-2}||_F^2 \\
&=& O(h^{-2r})\frac{\sigma^2}{n}O\left\{\frac{1}{\max(h, \lambdan^{1/2m})}\right\} \\
&=& O(h^{-2r})\frac{\sigma^2}{n}O\left\{\min(h^{-1}, \lambdan^{-1/2m})\right\}\\
&=& O(K^{2r})\frac{\sigma^2}{n}O\left\{\min(K, \lambdan^{-1/2m})\right\}\\
&=& O\left(\frac{K_e}{n}\right)
\end{eqnarray*}

Where in the above, we have used $||(I+\tilde{P})^{-2}||_F^2 = O\left\{\frac{1}{\max(h, \lambda^{1/2m})}\right\}$ from Lemma 5.2 of \cite{xiao2019}, $K\sim h^{-1}$, and $K_e = \min\left\{K^{2r+1}, K^{2r}\lambdan^{-1/2m}\right\}$

This completes the proof of the theorem.

\subsection{\textbf{Technical Lemmas}}

\begin{lemma}
\label{lm1}
    
Let $f\in\mathcal{C}^p$, then there exists $s_f\in \Sqt$ such that
$$||f^{(r)}-s_f^{(r)}|| = O(h^{q-r}) + o(h^{p-r})$$
for all $r = 0, 1, \dots, q - 2$ and $p \le q$.

Here, $b(x) = -\frac{f^{(q)}(x)h_i^q}{q!}B_q\left(\frac{x-t_i}{h_i}\right)$ where $B_q(.)$ is the $q^{th}$ Bernoulli polynomial defined as 
$B_0(x)=1$, and $B_k(x) = \displaystyle \int_0^x kB_{k-1}(x)dx + B_k$

and $B_k$ is chosen such that $\int_0^1B_k(x)dx=0$. \\
$B_k$ is known as the $k^{th}$ Bernoulli number \citep{barrow1978}. This Lemma also appears in \cite{xiao2019} where the general result in \cite{barrow1978} is adapted to prove the case where $p<q$.
\end{lemma}

\subsubsection*{Proof of Lemma \ref{lm1}}
We provide a proof for the case where $p=q$ and refer to Remark 3.1 of \cite{xiao2019} for the case where $p < q$. \cite{xiao2019} showed that when $p<q$, $||f^{(r)}-s_f^{(r)}|| = o(h^{p-r})$.\\
For $p=q$, first note that under Assumption \ref{assumption:3}, \cite{barrow1978} showed that 
$$\displaystyle \inf_{s(x)\in \Sqt}||\fp - s^{(r)}(x) + b^{*(r)}(x)||_{L_\infty} = o(h^{q-r})$$
This means, there exists an $s_f(x) \in \Sqt$ such that
$$||\fp - \sfp + b^{*(r)}(x)|| = o(h^{q-r})$$
where $b^*(x) = -\frac{f^{(q)}(t_i)h_i^q}{q!}B_q\left(\frac{x-t_i}{h_i}\right)$ and $b^{*(r)}$ is the $r^{th}$ derivative of $b^*(x)$.

With $p=q$, we have that $f\in\mathcal{C}^q[0, 1]$. Therefore, from Taylor's theorem, $f^{(q)}(x) = f^{(q)}(t_i) + o(1)$.
$$\implies b(x) = b^*(x) + o(h^q)$$
The derivative of the Bernoulli polynomial of order k is given by $\bm{B}'_k(x) = \bm{B}_{k-1}(x)$ \citep{barrow1978}, it therefore follows that
$$b^{(r)}(x) = b^{*(r)}(x) + o(h^{q-r})$$ 
for $r=0, 1, 2, \dots, q - 2$.
But $||b^*|| = O(h^q)$ by definition, giving $||b^{(r)}|| = O(h^{q-r})$.

Combining this with the case where $p<q$, we have that $||f^{(r)}-s_f^{(r)}|| = O(h^{q-r})+o(h^{p-r})$ for all $p\le q$.
\vspace{15pt}
\begin{lemma}
\label{lm2}
    
Given $\gp=\int_0^1 \bp\bpq q(x)dx$,
$$||\gp||_\infty = O(h^{-2r})$$
\end{lemma}

\subsubsection*{Proof of Lemma \ref{lm2}}
Note that $\bpq = \bqq D^{(r)}$
\begin{eqnarray*}
\therefore \gp &=& \int_0^1 \bqq D^{(r)}D^{T(r)}\bqqt q(x)dx \\
&=& O(h^{-2r})\int_0^1 \bqq\bqqt q(x)dx \\
&=& O(h^{-2r})\times q_{\max}\\
&=& O(h^{-2r})
\end{eqnarray*}

Where $q_{\max} = \displaystyle \max_{x\in [0, 1]} q(x) < \infty$. Also, note that B-spline bases are bounded by 1 $\forall x\in [0, 1]$.
\vspace{15pt}

\begin{lemma}
    \label{lm3}

Let $\gn = \bt\bm{B}/n$ where $\bm{B}=[B(x_1), B(x_2), \dots, B(x_n)]^T \in \mathbb{R}^{n\times K}$ is a matrix of basis functions with each $B(x)\in \mathbb{R}^{K}$ being a vector of basis functions of order $q$ at $x$.

Then
$$||\gni||_\infty = O(h^{-1})$$
\end{lemma}

\subsubsection*{Proof of Lemma \ref{lm3}}
This Lemma is adapted from \cite{zhou1998} and the key idea is to show that the elements of $\gni$ decay exponentially and of order $h^{-1}$. We provide the proof here for convenience. 

Let $\lambda_{\max} = \displaystyle \max_{ \sum_{i=1}^K a_i^2 = 1 } ||\gn \ab ||_2$ and $\lambda_{\min} = \displaystyle \min_{ \sum_{i=1}^K a_i^2 = 1 } ||\gn\ab||_2$ be the maximum and minimum eigen values of $\gn$.

Since $\gn$ is a band matrix, we use Theorem 2.2 of \cite{demko1977} by showing that the conditions of the theorem are satisfied.

First, note that
\begin{eqnarray*}
||\lambda^{-1}_{\max}\gn||_2 &=& \lambda^{-1}_{\max}||\gn||_2 \\
&=& \lambda^{-1}_{\max} \displaystyle \max_{ \sum_{i=1}^K z_i^2 = 1 } ||\gn\z||_2 \\
&\le& 1
\end{eqnarray*}

Where the $\max$ term in the second equality gives some eigen value which is at most $\lambda_{\max}^{-1}$.

Also, 
\begin{eqnarray*}
||\lambda_{\max}\gni||_2 &=& \frac{\lambda_{\max}}{\lambda_{\min}}||\lambda_{\min}\gni||_2\\
&\le& \frac{\lambda_{\max}}{\lambda_{\min}}
\end{eqnarray*}

Lemma 6.2 of \cite{zhou1998} provides bounds on the eigen values of $\gn$. In particular, for large $n$, there exist constants $c_1$ and $c_2$ such that 
$$c_1h/2 \le\lambda_{\min} \le \lambda_{\max} \le 2c_2h$$
Therefore by Theorem 2.2 of \cite{demko1977}, there exists constants $c > 0$ and $\gamma\in (0, 1)$ which depend only on $c_1$, $c_2$ and $q$ such that: 
\begin{equation}
\label{eqn:gni}
    |\lambda_{\max} g_{ij}| \le c\gamma^{|i-j|}
\end{equation}

where $g_{ij}$ is the $(i,j)$th element of $\gni$.

From equation \eqref{eqn:gni}, 
$$|g_{ij}| \le c\lambda_{\max}^{-1}\gamma^{|i-j|}\le 2(c/c_1)h^{-1}\gamma^{|i-j|}$$
This completes the proof of Lemma \ref{lm3}.

\begin{lemma}
\label{lm4}
Suppose $\gammab = \gni\bt\f/n$ and $P_m$ is the penalty matrix for the penalized spline estimator in \eqref{eqn:naive},

then
$$\gammabt P_m \gammab = O(1)$$
\end{lemma}

\subsubsection*{Proof of Lemma \ref{lm4}}
This Lemma is adapted from Lemma 8.4 of \cite{xiao2019} which puts a bound on the penalty matrix of the penalized spline estimator.
The proof follows closely the proof from \cite{xiao2019}.

Observe that 
\begin{eqnarray*}
\gni\bt\f/n &=& \gni\bt(\f-\sfb)/n + \gni\bt\sfb/n \\
&=& \gni\bt(\f-\sfb)/n + \betab \\
&=& \gni\alphab + \betab
\end{eqnarray*}
where $\betab = \gni\bt\sfb/n$ and $\alphab = \bt(\f-\sfb)/n$.
\begin{eqnarray}
\label{eqn:gamma2}
\left(\gammabt P_m\gammab \right)^{\frac{1}{2}} &\le& \left(\alphabt\gni P_m\gni\alphab\right)^{\frac{1}{2}} + \left(\betab^T P_m\betab\right)^{\frac{1}{2}}
\end{eqnarray}
since $P_m$ is positive semi-definite.
By Assumption, $\betab^TP_m\betab = O(1)$, therefore showing that the first term in \ref{eqn:gamma2} is O(1) completes the proof.

In the following, we use the following matrix relations.
Let $A\in \mathbb{R}^{m\times n}$, then
\begin{equation}
\label{mat:identity1}
    \frac{1}{\sqrt{n}}||A||_\infty\le||A||_2\le \sqrt{m}||A||_\infty
\end{equation}

Also, let $\pmh$ be a square symmetric matrix such that $P_m = \pmh\pmh$.

Observe that 
\begin{eqnarray}
\alphabt\gni P_m\gni\alphab &=& \left(\pmh\gni\alphab\right)^T\left(\pmh\gni\alphab\right)\\
&=& ||\pmh\gni\alphab||_2^2 \\
&\le& ||\alphab||_2^2||\pmh\gni||_2^2 \\
&\le& ||\alphab||_2^2||\pmh||_2^2||\gni||_2^2 \\
&\le& ||\alphab||_2^2||\pmh||_2^2K||\gni||_\infty^2 \\
&=& o(h^{2p+2})O(h^{1-2m})O(h^{-1})O(h^{-2}) \\
&=& o(h^{2p-2m}) \\
&=& o(1)
\end{eqnarray}
since $p\ge m$. Inequalities (12) and (13) are by Cauchy Schwartz inequality, we have used the matrix identity in \eqref{mat:identity1} in inequality (14). Also, We have used the result by \cite{agarwal1980} for $||\alphab||_2^2$ and the assumption that $||P_m||_2 = O(h^{1-2m})$. Finally, Lemma \ref{lm3} has been used in inequality (15) for $||\gni||_\infty$.

\subsection{Rates of Convergence for Naive Local Polynomial Estimators}
When estimating the $r^{th}$ derivative of the mean regression function with a local polynomial of degree $p$, several authors (\citealp[]{fan1996local, ruppert_wand_1994}) recommend using odd $p-r$. In this section, we lay out an argument that the naive bandwidth under- or over- smooths when $p$ and $p-r$ have different parities and that only even derivatives can be optimally estimated by the naive estimator. Table \ref{table:parities} below shows the four (4) potential parity combinations for $p$ and $p-r$. We show next that the naive estimator achieves the optimal rate of convergence when used to estimate $p-r$ only for cases I and IV (where $p$ and $p-r$ have same parity, equivalently, when $r$ is even).
\begin{table}[hbt!]
     \centering
     \begin{tabular}{|c|lcc|}
        \hline
        \multicolumn{4}{|r|}{$\bm{p-r}$} \\
        \hline
         ~&~& odd & even \\
         \multirow{2}{*}{$\bm{p}$} & odd & $I$ & $II$\\
         ~& even & $III$ & $IV$ \\
        \hline
        
    \end{tabular}
    \caption{Parity combinations of $p$ and $p-r$ when estimating the $r^{th}$ derivative of a mean regression function with $p^{th}$ degree local polynomial regression.}
    \label{table:parities}
\end{table}
Let $\mrh$ be a $p^{th}$-degree local polynomial estimate of the $r^{th}$ ($r\le p)$ derivative of the mean regression function, $m(x)$ at a point $x$ such that $m^{(p+1)}(\cdot)$ is continuous in a neighborhood of $x$. Let also $h$ be the bandwidth of $\mrh$ such that $h = o(n)$ and $nh \to \infty$, then we know from \cite{ruppert_wand_1994} that
$$\text{IMSE}\left(\mrh\right) = o\left(h^{2(p+1-r)}\right) + O\left(\frac{1}{nh^{2r+1}}\right)$$
for odd $p-r$ and
$$\text{IMSE}\left(\mrh\right) = o\left(h^{2(p+2-r)}\right) + O\left(\frac{1}{nh^{2r+1}}\right)$$
for even $p-r$.

Note that the naive estimator uses the optimal bandwidth when estimating $m(\cdot)$ itself, thus, when $r = 0$. In the above, IMSE is the integrated mean squared error.
First, we will derive the rates of convergence for the optimal bandwidth for the naive estimator $(r=0)$ for both the odd $p$ and even $p$ cases. We will then compare how these naive rates of convergence compare with the optimal bandwidths for estimating $p-v$ for both parity scenarios.

For odd $p$ (thus, $r=0$ and $p-r$ is odd),
$$\text{IMSE}\left(\mh\right) = o\left(h^{2(p+1)}\right) + O\left(\frac{1}{nh}\right)$$

To get the rate of convergence of the optimal bandwidth, we derive the $h$ that minimizes the IMSE (ignoring constants).

From:
\begin{eqnarray*}
2(p+1)h^{2p+1} - n^{-1}h^{-2} &=& 0 \\
2(p+1)h^{2p+1} &=& \frac{1}{nh^2} \\
h^{2p+3} &=& \frac{n^{-1}}{2(p+1)}
\end{eqnarray*}

$\therefore \hon = O\left(n^{-\frac{1}{2p+3}}\right)$. Here, we use $\hon$ for the optimal bandwidth for the naive estimator when $p$ is odd.

For even $p$ (thus, $r=0$ and $p-r$ is even),
$$\text{IMSE}\left(\mh\right) = o\left(h^{2(p+2)}\right) + O\left(\frac{1}{nh}\right)$$

From:
\begin{eqnarray*}
2(p+2)h^{2p+3} - n^{-1}h^{-2} &=& 0 \\
2(p+2)h^{2p+3} &=& \frac{1}{nh^2} \\
h^{2p+5} &=& \frac{n^{-1}}{2(p+2)}
\end{eqnarray*}

$\therefore \hen = O\left(n^{-\frac{1}{2p+5}}\right)$. $\hen$ is the optimal bandwidth for the naive estimator when $p$ is even.

We now analyse the achieved rates of convergence for estimating the $r^{th}$ derivative of the mean regression function, $m$ and how those rates compare with the naive estimator. We consider the four (4) cases in Table \ref{table:parities} above.

\textit{Case I}: \textbf{$p$ odd and $p-r$ odd (thus, $r$ is even).}
 
 $$\text{IMSE}\left(\mrh\right) = o\left(h^{2(p+1-r)}\right) + O\left(\frac{1}{nh^{2r+1}}\right)$$
 
 From 
 \begin{eqnarray*}
 2(p+1-r)h^{2p-2r+1} - (2r+1)n^{-1}h^{-2r-2} &=& 0\\
 2(p+1-r)h^{2p-2r+1} &=& \frac{2r+1}{nh^{2r+2}} \\
 h^{2p+3} &=& \frac{2r+1}{2(p+1-r)n}
 \end{eqnarray*}
 
 $\therefore h_{opt} = O\left(n^{-\frac{1}{2p+3}}\right)$, this is the same rate achieved by $\hon$. Therefore, the naive bandwidth achieves the same rate as the optimal bandwidth for estimating $p-r$ in this case.
 
 \textit{Case II}: \textbf{$p$ odd and $p-r$ even (thus, $r$ is odd).}
 
 $$\text{IMSE}\left(\mrh\right) = o\left(h^{2(p+2-r)}\right) + O\left(\frac{1}{nh^{2r+1}}\right)$$
By similar approach as in Case I, we get $h_{opt} = O\left(n^{-\frac{1}{2p+5}}\right)$, this rate is different from that achieved by the naive estimator $\hon$ for odd $p$. The consequence of using the naive bandwidth in this case is that, it shrinks faster than the optimal rate which may result in over-smoothing.
 
 \textit{Case III}: \textbf{$p$ even and $p-r$ odd (thus, $r$ is odd).}
 $$\text{IMSE}\left(\mrh\right) = o\left(h^{2(p+1-r)}\right) + O\left(\frac{1}{nh^{2r+1}}\right)$$
 
Again, similar to Cases I and II above, $h_{opt} = O\left(n^{-\frac{1}{2p+3}}\right)$, this rate is different from that achieved by the naive estimator $\hen$ for even $p$ which is $O\left(n^{-\frac{1}{2p+5}}\right)$. Unlike in case II, the consequence of using the naive bandwidth in this case is that, it shrinks at a slower rate than the optimal rate which may result in over-smoothing.
 
 \textit{Case IV}: \textbf{$p$ even and $p-r$ even (thus, $r$ is even).}
  $$\text{IMSE}\left(\mrh\right) = o\left(h^{2(p+2-r)}\right) + O\left(\frac{1}{nh^{2r+1}}\right)$$
 
 From 
 \begin{eqnarray*}
 2(p+2-r)h^{2p-2r+3} - (2r+1)n^{-1}h^{-2r-2} &=& 0\\
 2(p+2-r)h^{2p-2r+3} &=& \frac{2r+1}{nh^{2r+2}} \\
 h^{2p+5} &=& \frac{2r+1}{2(p+2-r)n}
 \end{eqnarray*}
 
 $\therefore h_{opt} = O\left(n^{-\frac{1}{2p+5}}\right)$, this is the same rate achieved by $\hen$ for even $p$. Therefore, the naive bandwidth achieves the same rate as the optimal bandwidth for estimating $p-r$ in this case. Thus, the naive estimator can only optimally estimate even-order derivatives for Local Polynomial Regression.